\documentclass{article}

\usepackage{times}
\usepackage{graphicx} % more modern
\usepackage{subfigure} 

\usepackage{natbib}

\usepackage{algorithm}
\usepackage{algorithmic}

\usepackage{hyperref}

\usepackage[accepted]{icml2017}
\usepackage{amsthm}
\usepackage{amsmath,amssymb}
\usepackage{subfigmat}

\theoremstyle{plain}
\newtheorem{thm}{Theorem}[section]
\newtheorem{lem}[thm]{Lemma}
\newtheorem{prop}[thm]{Proposition}
\newtheorem{cor}[thm]{Corollary}

\theoremstyle{definition}
\newtheorem{defn}{Definition}[section]

\theoremstyle{remark}
\newtheorem*{rem}{Remark}

\newtheorem*{pr}{Proof}
\newtheorem*{pr_main_thm_ns}{Proof of Theorem \ref{main_thm_ns}}

\newtheorem*{pr_main_thm_ns_warm}{Proof of Theorem \ref{main_thm_ns_warm}}
\newtheorem*{pr_main_cor_ns}{Proof of Corollary \ref{main_cor_ns}}
\newtheorem*{pr_main_lem_dasvrg}{Proof of Lemma \ref{main_lemma} for DASVRG}

\theoremstyle{assumption}
\newtheorem{assump}{Assumption}

\icmltitlerunning{DASVRDA for Regularized ERM}

\begin{document} 

\twocolumn[
\icmltitle{Doubly Accelerated Stochastic Variance Reduced Dual Averaging Method \\for Regularized Empirical Risk Minimization}



% you can specify symbols, otherwise they are numbered in order
\icmlsetsymbol{equal}{*}

\begin{icmlauthorlist}
\icmlauthor{Tomoya Murata}{gr}
\icmlauthor{Taiji Suzuki}{sc}
\end{icmlauthorlist}

\icmlaffiliation{gr}{NTT DATA Mathematical Systems Inc. Tokyo, Japan}
\icmlaffiliation{sc}{Graduate School of Information Science and Technology, The University of Tokyo}
%\icmlaffiliation{ed}{University of Edenborrow, Edenborrow, United Kingdom}

\icmlcorrespondingauthor{Tomoya Murata}{murata@msi.co.jp}
\icmlcorrespondingauthor{Taiji Suzuki}{taiji@mist.i.u-tokyo.ac.jp}

\icmlkeywords{convex optimization, empirical risk minimization, stochastic variance reduction, double acceleration, mini-batch, parallelization}

\vskip 0.3in
]



\printAffiliationsAndNotice{}  % leave blank if no need to mention equal contribution
%\footnotetext{hi}

\begin{abstract}
In this paper, we develop a new accelerated stochastic gradient method for efficiently solving the convex regularized empirical risk minimization problem in mini-batch settings. The use of mini-batches is becoming a golden standard in the machine learning community, because mini-batch settings stabilize the gradient estimate and can easily make good use of parallel computing. The core of our proposed method is the incorporation of our new ``double acceleration'' technique and variance reduction technique. We theoretically analyze our proposed method and show that our method much improves the mini-batch efficiencies of previous accelerated stochastic methods, and essentially only needs size $\sqrt{n}$ mini-batches for achieving the optimal iteration complexities for both non-strongly and strongly convex objectives, where $n$ is the training set size. Further, we show that even in non-mini-batch settings, our method achieves the best known convergence rate for both non-strongly and strongly convex objectives.
\end{abstract}

\section{Introduction}%%%%%%%%%%%%%%%%%%%%%%%%%%%%%%%%%%%%%%%%%%%%%
We consider a composite convex optimization problem associated with regularized empirical risk minimization, which often arises in machine learning. In particular, our goal is to minimize the sum of finite smooth convex functions and a relatively simple (possibly) non-differentiable convex function by using first order methods in mini-batch settings. The use of mini-batches is becoming a golden standard in the machine learning community, because it is generally more efficient to execute matrix-vector multiplications over a mini-batch than an equivalent amount of vector-vector ones each over a single instance; and more importantly, mini-batch settings can easily make good use of parallel computing. \par
Traditional and effective methods for solving the abovementioned problem are the ``proximal gradient'' (PG) method and ``accelerated proximal gradient'' (APG) method \cite{Nesterov2007gradient,beck2009fast,tseng2008accelerated}. These methods are well known to achieve linear convergence for strongly convex objectives. Particularly, APG achieves optimal iteration complexities for both non-strongly and strongly convex objectives. However, these methods need a per iteration cost of $O(nd)$, where $n$ denotes the number of components of the finite sum, and $d$ is the dimension of the solution space. In typical machine learning tasks, $n$ and $d$ correspond to the number of instances and features respectively, which can be very large. Then, the per iteration cost of these methods can be considerably high. \par
A popular alternative is the ``stochastic gradient descent'' (SGD) method \cite{singer2009efficient,hazan2007logarithmic,shalev2007logarithmic}. As the per iteration cost of SGD is only $O(d)$ in non-mini-batch settings, SGD is suitable for many machine learning tasks. However, SGD only achieves sublinear rates and is ultimately slower than PG and APG. \par

Recently, a number of stochastic gradient methods have been proposed; they use a variance reduction technique that utilizes the finite sum structure of the problem (``stochastic averaged gradient'' (SAG) method  \cite{roux2012stochastic,schmidt2013minimizing}, ``stochastic variance reduced gradient" (SVRG) method \cite{johnson2013accelerating,xiao2014proximal} and SAGA \cite{defazio2014saga}). Even though the per iteration costs of these methods are same as that of SGD, they achieve a linear convergence for strongly convex objectives. Consequently, these methods dramatically improve the total computational cost of PG. However, in size $b$ mini-batch settings, the rate is essentially $b$ times worse than in non-mini-batch settings. This means that there is little benefit in applying mini-batch scheme to these methods. \par
%More recently, several authors have proposed accelerated stochastic methods for the composite finite sum problem. \citet{shalev2013stochastic} has proposed Accelerated Stochastic Dual Coordinate Ascent (ASDCA) method, that can be seen as an accelerated inexact proximal point algorithm \cite{guler1992new,he2012accelerated,salzo2012inexact} where each subproblem is solved by (proximal) SDCA \cite{shalev2013stochastic}. Later, \citet{lin2015universal} has generalized this idea and shown that any linear convergent algorithms can be applied to this framework, and they call this generic scheme Universal Catalyst (UC). \citet{lin2014accelerated} have proposed Accelerated Proximal Coordinate Gradient (APCG) method and applied APCG to the dual problem. \citet{zhang2015stochastic} have proposed Stochastic Primal-Dual Coorinate (SPDC) based on a reformulation of the problem as a convex-concave saddle point problem. \citet{allen2016katyusha} has proposed Katyusha, that is the first primal-only accelerated stochastic variance reduction method, that is based on his new technique ``Katyusha momentum.'' 
More recently, several authors have proposed accelerated stochastic methods for the composite finite sum problem (``accelerated stochastic dual coordinate ascent'' (ASDCA) method \cite{shalev2013stochastic}, Universal Catalyst (UC) \cite{lin2015universal}, ``accelerated proximal coordinate gradient'' (APCG) method  \cite{lin2014accelerated}, ``stochastic primal-dual coordinate'' (SPDC) method \cite{zhang2015stochastic}, and Katyusha \cite{allen2016katyusha}). 
ASDCA (UC), APCG, SPDC and Katyusha essentially achieve the optimal total computational cost\footnotemark\footnotetext{More precisely, the rate of ASDCA  (UC) is with extra log-factors, and near but worse than the one of APCG, SPDC and Katyusha. This means that ASDCA (UC) cannot be optimal.} for strongly convex objectives\footnotemark\footnotetext{Katyusha also achieves  a near optimal total computational cost for non-strongly convex objectives.} in non-mini-batch settings. However, in size $b$ mini-batch settings, the rate is essentially $\sqrt{b}$ times worse than that in non-mini-batch settings, and these methods need size $O(n)$ mini-batches for achieving the optimal iteration complexity, which is essentially the same as APG. In addition, \citet{nitanda2014stochastic,nitanda2015accelerated} has proposed the ``accelerated mini-batch proximal stochastic variance reduced gradient'' (AccProxSVRG) method and its variant, the ``accelerated efficient mini-batch stochastic variance reduced gradient'' (AMSVRG) method. In non-mini-batch settings, AccProxSVRG only achieves the same rate as SVRG. However, in mini-batch settings, AccProxSVRG significantly improves the mini-batch efficiency of non-accelerated variance reduction methods, and surprisingly, AccProxSVRG essentially only needs size $O(\sqrt{\kappa})$ mini-batches for achieving the optimal iteration complexity for strongly convex objectives, where $\kappa$ is the condition number of the problem. However, the necessary size of mini-batches depends on the condition number and gradually increases when the condition number increases and ultimately matches with $O(n)$ for a large condition number. 
%Note that it is generally much more important to achieve the optimal iteration complextity for bad conditioned problems than well conditioned ones, because the former needs more overall computational cost than the latter. 

\begin{table*}[t]
\centering
\scalebox{0.67}[0.67]{
\hspace{-0.2cm}
\begin{tabular}{|c|c|c|c|c|c|c|}
\hline
& \multicolumn{3}{c|}{$\mu$-strongly convex} & \multicolumn{3}{c|}{Non-strongly convex} \\ \cline{2-7}
&Total computational cost &\multicolumn{2}{c|}{Necessary size of mini-batches}&Total computational cost &\multicolumn{2}{c|}{Necessary size of mini-batches} \\ \cline{3-4 } \cline{6-7} 
&in size $b$ mini-batch settings& $L/\mu \geq n$ & $L/\mu \leq n$&in size $b$ mini-batch settings& $L/\varepsilon \geq n\mathrm{log}^2(1/\varepsilon)$ & $L/\varepsilon \leq n\mathrm{log}^2(1/\varepsilon)$ \\ \hline
SVRG (SVRG$^{++}$) &$O\left( d\left(n + \frac{bL}{\mu}\right)\mathrm{log}\left(\frac{1}{\varepsilon}\right)\right)$&Unattainable&Unattainable&$O\left( d\left(n\mathrm{log}\left(\frac{1}{\varepsilon}\right) + \frac{bL}{\varepsilon}\right)\right)$& Unattainable &Unattainable\\ 
ASDCA (UC) &$\widetilde{O}\left(d\left(n+\sqrt{\frac{nbL}{\mu}}\right)\mathrm{log}\left(\frac{1}{\varepsilon}\right)\right)$ & Unattainable&Unattainable&$\widetilde{O}\left(d\left(\frac{n+\sqrt{nbL}}{\sqrt{\varepsilon}}\right)\right)$& Unattainable&Unattainable \\
APCG &$O\left(d\left(n+\sqrt{\frac{nbL}{\mu}}\right)\mathrm{log}\left(\frac{1}{\varepsilon}\right)\right)$&$O(n)$&$O(n)$ &No direct analysis&Unattainable&Unattainable\\ 
SPDC&$O\left(d\left(n+\sqrt{\frac{nbL}{\mu}}\right)\mathrm{log}\left(\frac{1}{\varepsilon}\right)\right)$&$O(n)$&$O(n)$&No direct analysis&Unattainable&Unattainable\\ 
Katyusha &$O\left(d\left(n+\sqrt{\frac{nbL}{\mu}}\right)\mathrm{log}\left(\frac{1}{\varepsilon}\right)\right)$&$O(n)$&$O(n)$&$O\left( d\left(n\mathrm{log}\left(\frac{1}{\varepsilon}\right)+\sqrt{\frac{nbL}{\varepsilon}}\right)\right)$&$O(n)$&$O(n)$\\
AccProxSVRG &$O\left(d\left( n + \left(\frac{n-b}{n-1}\right)\frac{L}{\mu} + b\sqrt{\frac{L}{\mu}}\right)\mathrm{log}\left(\frac{1}{\varepsilon}\right)\right)$&$O\left(\sqrt{\frac{L}{\mu}}\right)$ &$O\left(n\sqrt{\frac{\mu}{L}}\right)$ &No direct analysis&Unattainable&Unattainable\\ 
%UC $+$ SVRG &$O\left(\left(n+\sqrt{\frac{nbL}{\mu}}\right)\mathrm{log}\biggl(\frac{1 + \frac{bL}{n\mu}}{\varepsilon}\biggr) \mathrm{log}\left(1 + \frac{bL}{n\mu}\right) \right)$&Near opt. &$O\left(\frac{n + \sqrt{nbL}}{\sqrt{\varepsilon}} \mathrm{log}\left(1+\frac{bL}{n}\Bigl(\frac{1+\frac{bL}{n}}{\varepsilon}\Bigr)^{4+\eta}\right) \right)$&Near opt. \\
%UC $+$ AccProxSVRG\footnotemark[3]&$O\left(\left(n+\sqrt{\frac{nL}{\mu}}+b\sqrt{\frac{L}{\mu}}\right)\mathrm{log}\biggl(\frac{ 1 + \frac{L}{n\mu}}{\varepsilon}\biggr) \mathrm{log}\left( 1 + \frac{L}{n\mu}\right) \right)$&Near opt. &$O\left(\frac{n + \sqrt{nL} + b\sqrt{L}}{\sqrt{\varepsilon}}\mathrm{log}\left( 1 + \frac{L}{n}\Bigl( \frac{1+\frac{L}{n}}{\varepsilon}\Bigr)^{4+\eta} + \frac{b^2L}{n^2}\Bigl( \frac{1+\frac{b^2L}{n^2}}{\varepsilon}\Bigr)^{4+\eta}\right)\right)$&Near opt. \\
{\bfseries{DASVRDA}} &{\color{red}$O\left(d\left(n+\sqrt{\frac{nL}{\mu}} + b\sqrt{\frac{L}{\mu}}\right)\mathrm{log}\left(\frac{1}{\varepsilon}\right)\right)$}&{\color{red}$O\left(\sqrt{n}\right)$} &$O\left(n\sqrt{\frac{\mu}{L}}\right)$
&{\color{red}$O\left(d\left(n\mathrm{log}\frac{1}{\varepsilon}+\sqrt{\frac{nL}{\varepsilon}} + b\sqrt{\frac{L}{\varepsilon}}\right)\right)$}&{\color{red}$O\left(\sqrt{n}\right)$}& $\widetilde O\left(n\sqrt{\frac{\varepsilon}{L}}\right)$\\ \hline
\end{tabular}}
\caption{Comparisons of our method with SVRG (SVRG$^{++}$ \cite{AllenYang2016}), ASDCA (UC), APCG, SPDC, Katyusha and AccProxSVRG. $n$ is the number of components of the finite sum, $d$ is the dimension of the solution space, $b$ is the mini-batch size, $L$ is the smoothness parameter of the finite sum, $\mu$ is the strong convexity parameter of objectives (see Def. \ref{L-smooth_def} and Def. \ref{sc_def} for their definitions), and $\varepsilon$ is accuracy.  ``Necessary size of mini-batches '' indicates the order of the necessary size of mini-batches for achieving optimal iteration complexities $O(\sqrt{L/\mu}\mathrm{log}(1/\varepsilon))$ and $O(\sqrt{L/\varepsilon})$ for strongly and non-strongly convex objectives, respectively. We regard one computation of a full gradient as $n/b$ iterations in size $b$ mini-batch settings, for a fair comparison. ``Unattainable'' implies that the algorithm cannot achieve the optimal iteration complexity even if it uses size $n$ mini-batches. $\widetilde{O}$ hides extra log-factors. The results marked in red denote the main contributions of this paper. 
}
\label{table}
\end{table*} 
\subsection*{Main contribution}
We propose a new accelerated stochastic variance reduction method 
that achieves better convergence than existing methods do, and 
it particularly utilizes mini-batch settings well; 
it is called the ``doubly accelerated stochastic variance reduced dual averaging'' (DASVRDA) method. 
%We describe the main features of our proposed method below and 
Our method significantly improves the mini-batch efficiencies of state-of-the-art methods, and our method essentially only needs size $O(\sqrt{n})$ mini-batches\footnotemark\footnotetext{Actually, when  $L/\varepsilon \leq n\mathrm{log}^2n$ and $L/\mu \leq n$, our method needs size $\widetilde{O}(n\sqrt{\varepsilon/L})$ and $O(n\sqrt{\mu/L})$ mini-batches, respectively, which are larger than $O(\sqrt{n})$, but smaller than $O(n)$. Achieving the optimal iteration complexity for solving high accuracy and bad conditioned problems is much more important than doing ones with low accuracy and well-conditioned ones, because the former needs more overall computational cost than the latter. } for achieving the optimal iteration complexities\footnotemark\footnotetext{We refer to ``optimal iteration complexity'' as the iteration complexity of deterministic Nesterov's acceleration method (\citet{nesterov2013introductory})} for both non-strongly and strongly convex objectives. 
We list the comparisons of our method with several preceding methods in Table \ref{table}. 

\section{Preliminary}%%%%%%%%%%%%%%%%%%%%%%%%%%%%%%%%%%%%%%%%%%%%%%
In this section, we provide several notations and definitions used in this paper. Then, we make assumptions for our analysis. \par
We use $\| \cdot \|$ to denote the Euclidean $L_2$ norm $\| \cdot \|_2$: $\|x\| = \|x\|_2 = \sqrt{\sum_{i}x_i^2}$. For natural number $n$, $[n]$ denotes set $\{1, \ldots, n\}$. 
\begin{defn}\label{L-smooth_def}
We say that a function $f: \mathbb{R}^d \to \mathbb{R}$ is  $L$-smooth ($L > 0$) if $f$ is differentiable and satisfies
\begin{equation}
\|\nabla f(x)-\nabla f(y)\| \leq L\|x-y\|\ (\forall x, y \in \mathbb{R}^d). \label{L-smooth_property}
\end{equation}
\end{defn}
If $f$ is convex, (\ref{L-smooth_property}) is equivalent to the following:  
(see \citet{nesterov2013introductory}): 
\begin{align}
&f(x) +  \langle \nabla f(x), y-x\rangle + \frac{1}{2L}\|\nabla f(x) - \nabla f(y)\|^2 \leq f(y) \label{L-smooth_ineq} \\
&(\forall x, y \in \mathbb{R}^d). \notag
\end{align}
\begin{defn}\label{sc_def}
A convex function $f: \mathbb{R}^d \to \mathbb{R}$ is called $\mu$-strongly convex ($\mu \geq 0$) if $f$ satisfies
\begin{align}
&\frac{\mu}{2}\|x-y\|^2 + \langle \xi, y - x \rangle + f(x) \leq f(y), \label{strongly_convex_def} \notag \\
&(\forall x, y \in \mathbb{R}^d, \forall \xi \in \partial f(x)) \notag 
\end{align}
where $\partial f(x)$ denotes the set of the subgradients of $f$ at $x$. 
\end{defn}
Note that if $f$ is $\mu$-strongly convex, then $f$ has the unique minimizer. 
\begin{defn}\label{osc_def}
We say that a function $f: \mathbb{R}^d \to \mathbb{R}$ is $\mu$-optimally strongly convex ($\mu \geq 0$) if $f$ has a minimizer and satisfies 
\begin{align}
&\frac{\mu}{2}\| x - x_* \|^2 \leq f(x) - f(x_*) \notag \\
&(\forall x \in \mathbb{R}^d, \forall x_* \in \mathrm{argmin}_{x \in \mathbb{R}^d} f(x)). \notag 
\end{align}
\end{defn}
If a function $f$ is $\mu$-strongly convex, then $f$ is clearly $\mu$-optimally strongly convex. 
\begin{defn}
We say that a convex function $f: \mathbb{R}^d \to \mathbb{R}$ is relatively simple if computing the proximal mapping of $f$ at $y$, 
$$\mathrm{prox}_{f}(y)=\underset{x \in \mathbb{R}^d}{\mathrm{argmin}\ } \left\{ \frac{1}{2}\|x-y\|^{2} + f(x) \right\},$$
takes at most $O(d)$ computational cost, for any $y \in \mathbb{R}^d$.   
\end{defn}
As $f$ is convex, function $(1/2)\|x-y\|^{2} + f(x)$ is $1$-strongly convex, and function $\mathrm{prox}_{f}$ is well-defined. \par
Now, we formally describe the problem to be considered in this paper and the assumptions for our theory. In this paper, we consider the following composite convex minimization problem:
\begin{equation}\label{problem}
\underset{x \in \mathbb{R}^d}{\mathrm{min}}\ \ \{P(x)\overset{\mathrm{def}}{=}F(x)+R(x)\}, 
\end{equation}
where  $F(x)=\frac{1}{n}\sum_{i=1}^n f_{i}(x)$. Here each $f_{i}:\mathbb{R}^d \to \mathbb{R}$ is a $L_i$-smooth convex function and $R:\mathbb{R}^d \to \mathbb{R}$ is a relatively simple and (possibly) non-differentiable convex function. Problems of this form
 often arise in machine learning and fall under regularized empirical risk minimization (ERM). In ERM problems, we are given $n$ training examples $\{ (a_i, b_i) \}_{i=1}^n$, where each $a_i \in \mathbb{R}^d$ is the feature vector of example $i$, and each $b_i \in \mathbb{R}$ is the label of example $i$. The following regression and classification problems are typical examples of ERM on our setting: 
\begin{itemize}
\item Lasso: $f_i(x) = \frac{1}{2}(a_i^{\top}x - b_i)^2$ and $R(x) = \lambda \|x\|_1$. 
\item Ridge logistic regression: $f_i(x) = \mathrm{log}(1+\mathrm{exp}(-b_ia_i^{\top}x))$ and $R(x) = \frac{\lambda}{2}\|x\|_2^2$. 
\item Support vector machines: $f_i(x) = \bar h_i^{\nu}(a_i^{\top}x)$ and $R(x) = \frac{\lambda}{2}\|x\|_2^2$.
\end{itemize}
Here, $\bar h_i^{\nu}$ is a smooth variant of hinge loss (for the definition, for example, see \citet{shalev2013stochastic}). \par
%$h_i(u) \overset{\mathrm{def}}{=}\mathrm{max}\{0, 1-b_i u\}$ and defined by 
%$$\bar h_i^{\nu} (u) \overset{\mathrm{def}}{=} 
%\begin{cases} 1 - b_i u - \frac{1}{2\nu} \ \ \ \ (b_i u < 1 - \frac{1}{\nu}) \\
%\frac{\nu}{2}(1 - b_i u)^2 \ \ \ \ \hspace{0.4em}(1-\frac{1}{\nu} \leq b_i u \leq 1) \\
%0 \ \ \ \ \hspace{4.9em}(b_i u > 1)
%\end{cases}
%,$$ where $\nu$ is the smoothing parameter and we can easily see that $\bar h_i^{\nu}$ is $\nu$-smooth
%. \par
We make the following assumptions for our analysis: 
\begin{assump}\label{optimal_solution_assump}
There exists a minimizer $x_*$ of (\ref{problem}).
\end{assump}
\begin{assump}\label{L-smooth_assump}
Each $f_i$ is convex and $L_i$-smooth.
\end{assump}
The above examples satisfy this assumption with $L_i = \|a_i\|_2^2$ for the squared loss, $L_i = \|a_i\|_2^2/4$ for the logistic loss, and $L_i = \nu \|a_i\|_2^2$ for the smoothed hinge loss. 
%We denote $\bar L$ as $\frac{1}{n}\sum_{i=1}^n L_i$. Moreover, we define the probability distribution $Q$ on the set $\{ 1, 2, \ldots, n \}$ by $Q = \{q_i\}_{i=1}^n \overset{\mathrm{def}}{=} \left\{ \frac{L_i}{n\bar L} \right\}_{i=1}^n$. 
\begin{assump}\label{proximal_assump}
Regularization function $R$ is convex and is relatively simple. 
\end{assump}
For example, Elastic Net regularizer $R(x) = \lambda_1 \|x\|_1 + (\lambda_2/2) \|x\|_2^2$ ($\lambda_1, \lambda_2 \geq 0$) satisfies this assumption. Indeed, we can analytically compute the proximal mapping of $R$ by $\mathrm{prox}_{R}(z)=((1/(1 + \lambda_2)) \mathrm{sign}(z_j)\mathrm{max} \{|z_j|-\lambda_1, 0\})_{j=1}^d$, and this costs only $O(d)$. \par
We always consider Assumption \ref{optimal_solution_assump}, \ref{L-smooth_assump} and \ref{proximal_assump} in this paper. 
\begin{assump}\label{osc_assump}
There exists $\mu > 0$ such that objective function $P$ is $\mu$-optimally strongly convex.
\end{assump}
If each $f_i$ is convex, then for Elastic Net regularization function $R(x) = \lambda_1 \|x\|_1 + (\lambda_2/2) \|x\|_2^2$ ($\lambda_1\geq 0$ and $\lambda_2 > 0$), Assumption \ref{osc_assump} holds with $\mu = \lambda_2$. \par
%Note that Assumption \ref{osc_assump} is weaker than Assumption \ref{sc_assump}. 
We further consider Assumption \ref{osc_assump} when we deal with strongly convex objectives.
\section{Our Approach: Double Acceleration}\label{approach_sec}%%%%%%%%%%%%%%%%%%%%%%%%%%%%%%%%%%%%%%%%%%%
In this section, we provide high-level ideas of our main contribution called ``double acceleration.'' \par %Here, we focus on non-strongly convex objectives. \par
\setlength\textfloatsep{0.4cm}
\begin{algorithm}[t]
\caption{PG $(\widetilde{x}_{0}, \eta, S$)}
\label{pg_ns}
\begin{algorithmic}
%[1]
\REQUIRE $\widetilde{x}_{0} \in \mathbb{R}^d$, $\eta>0$, $S \in  \mathbb{N}$. 
\FOR {$s=1$ to $S$}
\STATE $\widetilde{x}_s =  \text{One Stage PG}(\widetilde{x}_{s-1}, \eta). $
\ENDFOR
\ENSURE $\frac{1}{S}\sum_{s=1}^S\widetilde{x}_s$. 
\end{algorithmic}
\end{algorithm}
\setlength\floatsep{0.3cm}
\begin{algorithm}[t]
\caption{One Stage PG $(\widetilde{x}, \eta$)}
\label{one_Stage_pg_ns}
\begin{algorithmic}
%[1]
\REQUIRE $\widetilde{x} \in \mathbb{R}^d$, $ \eta > 0$. 
\STATE $\widetilde{x}^+ = \mathrm{prox}_{\eta R}(\widetilde{x} - \eta \nabla F(\widetilde{x}))$. 
\ENSURE $\widetilde{x}^+$. 
\end{algorithmic}
\end{algorithm}
\begin{algorithm}[t]
\caption{SVRG $(\widetilde{x}_{0}, \eta, m, S$)}
\label{svrg_ns}
\begin{algorithmic}
%[1]
\REQUIRE $\widetilde{x}_{0} \in \mathbb{R}^d$, $\eta>0$, $m, S \in  \mathbb{N}$. 
\FOR {$s=1$ to $S$}
\STATE $\widetilde{x}_s =  \text{One Stage SVRG}(\widetilde{x}_{s-1}, \eta, m). $
\ENDFOR
\ENSURE $\frac{1}{S}\sum_{s=1}^S\widetilde{x}_s$. 
\end{algorithmic}
\end{algorithm}
\begin{algorithm}[t]
\caption{One Stage SVRG $(\widetilde{x}, \eta, m$)}
\label{one_Stage_svrg_ns}
\begin{algorithmic}
%[1]
\REQUIRE $\widetilde{x} \in \mathbb{R}^d$, $ \eta > 0$, $m \in \mathbb{N}$. 
\STATE $x_0 = \widetilde{x}$. 
\FOR {$k=1$ to $m$}
\STATE Pick $i_k \in \{1, \ldots, n\}$ randomly.
\STATE $g_k = \nabla f_{i_k}(x_{k-1}) - \nabla f_{i_k}(\widetilde{x}) + \nabla F(\widetilde{x})$. \label{variance_reduction_svrg_ns}
\STATE $x_k = \mathrm{prox}_{\eta R}(x_{k-1} - \eta g_k)$. 
\ENDFOR
\ENSURE $\frac{1}{n}\sum_{k=1}^n x_k$. 
\end{algorithmic}
\end{algorithm}
First, we consider deterministic PG (Algorithm \ref{pg_ns}) and (non-mini-batch) SVRG (Algorithm \ref{svrg_ns}). PG is an extension of the steepest descent to proximal settings. SVRG is a stochastic gradient method using the variance reduction technique, which utilizes the finite sum structure of the problem, and it achieves a faster convergence rate than PG does. 
%Indeed, we can show that SVRG (Algorithm \ref{svrg_ns}) achieves a total computational cost of $O(d(n+L)/\varepsilon)$ in serial non-mini-batch settings if we choose $\eta = O(1/L)$ and $m = O(n)$, which is better than the one of PG (Since $m = O(n)$, we call Algorithm \ref{one_Stage_svrg_ns} ``One Stage'' SVRG). 
As SVRG (Algorithm \ref{svrg_ns}) matches with PG (Algorithm \ref{pg_ns}) when the number of inner iterations is $m=1$, SVRG can be seen as a generalization of PG. 
%Now, we briefly explain the variance reduction technique, that is essential for SVRG. 
The key element of SVRG is employing a simple but powerful technique called the variance reduction technique for gradient estimate.
%One Stage SVRG uses 
The variance reduction of the gradient is realized by setting $g_k = \nabla f_{i_k}(x_{k-1}) - \nabla f_{i_k}(\widetilde{x}) + \nabla F(\widetilde{x})$ rather than vanilla stochastic gradient $\nabla f_{i_k}(x_{k-1})$. Generally, stochastic gradient $\nabla f_{i_k}(x_{k-1})$ is an unbiased estimator of $\nabla F(x_{k-1})$, but it may have high variance. In contrast, $g_k$ is also unbiased, and one can show that its variance is ``reduced''; that is, the variance converges to zero as $x_{k-1}$ and $\widetilde{x}$ to $x_*$. \par

Next, we explain to the method of accelerating SVRG and obtaining an even faster convergence rate based on our new but quite natural idea ``outer acceleration.'' 
First, we would like to remind you that the procedure of deterministic APG is given as described in Algorithm \ref{apg_ns}. APG uses the famous ``momentum'' scheme and achieves the optimal iteration complexity. Our natural idea is replacing One Stage PG in Algorithm \ref{apg_ns} with One Stage SVRG. With slight modifications, we can show that this algorithm improves the rates of PG, SVRG and APG, and is optimal. We call this new algorithm outerly accelerated SVRG (Note that the algorithm matches with APG when $m=1$ and thus, can be seen as a direct generalization of APG). However, this algorithm has poor mini-batch efficiency, because in size $b$ mini-batch settings, the rate of this algorithm is essentially $\sqrt{b}$ times worse than that of non-mini-batch settings. State-of-the-art methods APCG, SPDC, and Katyusha also suffer from the same problem in the mini-batch setting. \par
\begin{algorithm}[t]
\caption{APG $(\widetilde{x}_{0}, \eta, S$)}
\label{apg_ns}
\begin{algorithmic}
%[1]
\REQUIRE $\widetilde{x}_{0} \in \mathbb{R}^d$, $\eta>0$, $S \in  \mathbb{N}$. 
\STATE $\widetilde{x}_{-1} = \widetilde{x}_0$. 
\STATE $\widetilde{\theta}_0 = 0$. 
\FOR {$s=1$ to $S$}
\STATE $\widetilde{\theta}_s = \frac{s+1}{2}$.
\STATE $\widetilde{y}_s = \widetilde{x}_{s-1} + \frac{\widetilde{\theta}_{s-1}-1}{\widetilde{\theta}_s}(\widetilde{x}_{s-1} - \widetilde{x}_{s-2})$. 
\STATE $\widetilde{x}_s =  \text{One Stage PG}(\widetilde{y}_{s}, \eta)$. \label{replace}
\ENDFOR
\ENSURE $x_S$. 
\end{algorithmic}
\end{algorithm}
Now, we illustrate that for improving the mini-batch efficiency, using the ``inner acceleration'' technique is beneficial. Nitanda \cite{nitanda2014stochastic} has proposed AccProxSVRG in mini-batch settings. He applied the momentum scheme to One Stage SVRG, and we call this technique ``inner'' acceleration.
He showed that the inner acceleration could significantly improve the mini-batch efficiency of vanilla SVRG. 
This fact indicates that inner acceleration is essential to fully utilize the mini-batch settings. 
However, AccProxSVRG is not a truly accelerated method, because in non-mini-batch settings, the rate of AccProxSVRG is same as that of vanilla SVRG. \par
In this way, we arrive at our main high-level idea called ``double'' acceleration, which involves applying momentum scheme to both outer and inner algorithms. This enables us not only to lead to the optimal total computational cost in non-mini-batch settings, but also to improving the mini-batch efficiency of vanilla acceleration methods. \par
We have considered SVRG and its accelerations so far; however, we actually adopt stochastic variance reduced dual averaging (SVRDA) rather than SVRG itself, because we can construct lazy update rules of (innerly) accelerated SVRDA for sparse data (see Section \ref{lazy_update_sec}). In Section \ref{dasvrg_sec} of supplementary material, we briefly discuss a SVRG version of our proposed method and provide its convergence analysis.  

\section{Algorithm Description}%%%%%%%%%%%%%%%%%%%%%%%%%%%%%%%%%%%%%%%%%%%%
\label{algorithm_sec}
In this section, we describe the concrete procedure of the proposed algorithm in detail. \par
\subsection*{DASVRDA for non-strongly convex objectives}
We provide details of the doubly accelerated SVRDA (DASVRDA) method for non-strongly convex objectives in Algorithm \ref{dasvrda_ns}. Our momentum step is slightly different from that of vanilla deterministic accelerated methods: we not only add momentum term $((\widetilde{\theta}_{s-1}-1)/\widetilde{\theta}_s)(\widetilde{x}_{s-1} - \widetilde{x}_{s-2})$ to the current solution $\widetilde{x}_{s-1}$ but also add term $(\widetilde{\theta}_{s-1}/\widetilde{\theta}_s)(\widetilde{z}_{s-1} - \widetilde{x}_{s-1})$, where $\widetilde{z}_{s-1}$ is the current more ``aggressively'' updated solution rather than $\widetilde{x}_{s-1}$; thus, this term also can be interpreted as momentum\footnotemark\footnotetext{This form also arises in Monotone APG \cite{li2015accelerated}. In Algorithm \ref{one_Stage_accsvrda}, $\widetilde{x}=x_m$ can be rewritten as $(2/(m(m+1))\sum_{k=1}^m kz_k$, which is a weighted average of $z_k$; thus, we can say that $\widetilde{z}$ is updated more ``aggressively'' than $\widetilde{x}$. 
For the outerly accelerated SVRG (that is a combination of Algorithm \ref{dasvrda_ns} with vanilla SVRG, see section \ref{approach_sec}), $\widetilde{z}$ and $\widetilde{x}$ correspond to $x_m$ and $(1/m)\sum_{k=1}^m x_k$ in \cite{xiao2014proximal}, respectively. Thus, we can also see that $\widetilde{z}$ is updated more ``aggressively'' than $\widetilde{x}$. 
}. Then, we feed $\widetilde{y}_s$ to One Stage Accelerated SVRDA (Algorithm \ref{one_Stage_accsvrda}) as an initial point. Note that Algorithm \ref{dasvrda_ns} can be seen as a direct generalization of APG, because if we set $m=1$, One Stage Accelerated SVRDA is essentially the same as one iteration PG with initial point $\widetilde{y}_s$; then, we can see that $\widetilde{z}_s = \widetilde{x}_s$, and Algorithm \ref{dasvrda_ns} essentially matches with deterministic APG. 
Next, we move to One Stage Accelerated SVRDA (Algorithm \ref{one_Stage_accsvrda}). Algorithm \ref{one_Stage_accsvrda} is essentially a combination of the ``accelerated regularized dual averaging'' (AccSDA) method \cite{xiao2009dual} with the variance reduction technique of SVRG. It updates $z_k$ by using the weighted average of all past variance reduced gradients $\bar g_k$ instead of only using a single variance reduced gradient $g_k$. Note that for constructing variance reduced gradient $g_k$, we use the full gradient of $F$ at $\widetilde{x}_{s-1}$ rather than the initial point $\widetilde{y}_s$. \par
%Line \ref{momentum_to} can be interpreted as the momentum step, because we can rewrite line \ref{momentum_to} as $y_k = x_{k-1} + ((\theta_{k-1}-1)/\theta_k)(x_{k-1} - x_{k-2})$ by using the relation in line \ref{convex_comb}.  
%In line \ref{orda_line}, we set $x_k$ as a convex combination of $x_{k-1}$ and $z_k$. 
%Optimal Regularized Dual Ageraging (ORDA) method \cite{chen2012optimal} sets $x_k$ as the output of a (proximal) SGD step at $y_{k}$ rather than take the convex combination. This step generates a sparse output solution $x_m$ for sparsity inducing regularizer, whereas taking the convex combination may make the output solution non-sparse. We can apply this scheme to our algorithm, however it seems to be difficult to construct lazy update rules of this algorithm, whereas using the convex combination enables us to construct ones. Thus, we adopt the convex combination scheme for our algorithm. 
\begin{rem}
In Algorithm \ref{one_Stage_accsvrda}, we pick $b$ indexes according to i.i.d. non-uniform distribution $Q = \{ q_i\} = \left\{ \frac{L_i}{n\bar L} \right\}$. Instead of that, we can pick $b$ indexes, such that each index $i_k^{\ell}$ is uniformly picked from $B^{\ell}$, where $\{ B^{\ell}\}_{\ell=1}^b$ is the predefined disjoint partition of $[n]$ with size $|B^{\ell}| = n/b$. If we adopt this scheme, when we parallelize the algorithm using $b$ machines, each machine only needs to store the corresponding partition of the data set rather than the whole dataset, and this can reduce communication cost and memory cost. In this setting, the convergence analysis in Section \ref{analysis_sec} can be easily revised by simply replacing $\bar L$ with $L_{\mathrm{max}} = \mathrm{max}_{i \in [n]}L_i$.
\end{rem}

\subsection*{DASVRDA for strongly convex objectives}
Algorithm \ref{dasvrda_sc} is our proposed method for strongly convex objectives. Instead of directly accelerating the algorithm using a constant momentum rate, we restart Algorithm \ref{dasvrda_ns}. 
%, that means starting Algorithm \ref{dasvrda_ns} again, taking the current solution as the new initial point
Restarting scheme has several advantages both theoretically and practically. First, the restarting scheme only requires the optimal strong convexity of the objective (Def. \ref{osc_def}) instead of the ordinary strong convexity (Def. \ref{sc_def}). Whereas, non-restarting accelerated algorithms essentially require the ordinary strong convexity of the objective. Second, for restarting algorithms, we can utilize adaptive restart schemes \cite{o2015adaptive}. The adaptive restart schemes have been originally proposed for deterministic cases. The schemes are heuristic but quite effective empirically. The most fascinating property of these schemes is that we need not prespecify the strong convexity parameter $\mu$, and the algorithms adaptively determine the restart timings.  \citet{o2015adaptive} have proposed two heuristic adaptive restart schemes: the function scheme and gradient scheme. We can easily apply these ideas to our method, because our method is a direct generalization of the deterministic APG. For the function scheme, we restart Algorithm \ref{dasvrda_ns} if $P(\widetilde{x}_s) > P(\widetilde{x}_{s-1})$. For the gradient scheme, we restart the algorithm if $(\widetilde{y}_s - \widetilde{x}_s)^{\top}(\widetilde{y}_{s+1} - \widetilde{x}_{s}) > 0$. Here $\widetilde{y}_s - \widetilde{x}_s$ can be interpreted as a ``one stage'' gradient mapping of $P$ at $\widetilde{y}_s$. As $\widetilde{y}_{s+1} - \widetilde{x}_{s}$ is the momentum, this scheme can be interpreted such that we restart whenever the momentum and negative one Stage gradient mapping form an obtuse angle (this means that the momentum direction seems to be ``bad''). 
%It is one of the advantages of our outer acceleration technique that one can naturally apply the gradient scheme to the algorithm. 
We numerically demonstrate the effectiveness of these schemes in Section \ref{experiment_sec}. \par

\subsection*{DASVRDA$^{\mathrm{ns}}$ with warm start}
Algorithms \ref{warm_start} is a combination of DASVRDA$^{\mathrm{ns}}$ with warm start scheme. At the warm start phase, we repeatedly run One Stage AccSVRDA and increment the number of its inner iterations $m_u$ exponentially until $m_u \propto n/b$. After that, we run vanilla DASVRDA$^{\mathrm{ns}}$. We can show that this algorithm gives a faster rate than vanilla DASVRDA$^{\mathrm{ns}}$. 
\begin{rem}
For DASVRDA$^{\mathrm{sc}}$, the warm start scheme for DASVRDA$^{\mathrm{ns}}$ is not needed because the theoretical rate is identical to the one without warm start. 
\end{rem}
\subsection*{Parameter tunings}
For DASVRDA$^{\mathrm{ns}}$, only learning rate $\eta$ needs to be tuned, because we can theoretically obtain the optimal choice of $\gamma$, and we can naturally use $m = n/b$ as a default epoch length (see Section \ref{analysis_sec}). For DASVRDA$^{\mathrm{sc}}$, both learning rate $\eta$ and fixed restart interval $S$ need to be tuned. 

\begin{algorithm}[t]
\caption{DASVRDA$^{\mathrm{ns}}(\widetilde{x}_{0}, \widetilde{z}_0, \gamma, \{ L_i \}_{i=1}^n , m, b, S$)}
\label{dasvrda_ns}
\begin{algorithmic}
%[1]
\REQUIRE $\widetilde{x}_{0}, \widetilde{z}_0 \in \mathbb{R}^d$, $ \gamma > 1$, $\{ L_i > 0 \}_{i=1}^n$, $m \in \mathbb{N}$, $S \in  \mathbb{N}$, $b \in [n]$. 
%\STATE $\eta = \frac{1}{\left(1 + \frac{\gamma(m+1)}{b}\right)\hat L}$. 
\STATE $\widetilde{x}_{-1} = \widetilde{z}_0$, $\widetilde{\theta}_0 = 1-\frac{1}{\gamma}$. 
\STATE $\bar L = \frac{1}{n}\sum_{i=1}^n L_i$.  
\STATE $Q = \{ q_i\} = \left\{ \frac{L_i}{n\bar L} \right\}$. 
\STATE $\eta = \frac{1}{\left(1 + \frac{\gamma (m+1)}{b}\right)\bar L}$. 
\FOR {$s=1$ to $S$}
\STATE $\widetilde{\theta}_s = \left(1 - \frac{1}{\gamma}\right)\frac{s+2}{2}. $
\STATE $\widetilde{y}_{s} =  \widetilde{x}_{s-1} + \frac{\widetilde{\theta}_{s-1}-1}{\widetilde{\theta}_s}(\widetilde{x}_{s-1} - \widetilde{x}_{s-2}) + \frac{\widetilde{\theta}_{s-1}}{\widetilde{\theta}_s}(\widetilde{z}_{s-1} - \widetilde{x}_{s-1}). $
\STATE $(\widetilde{x}_s, \widetilde{z}_s) =  \text{One Stage AccSVRDA}(\widetilde{y}_s, \widetilde{x}_{s-1}, \eta, m, $ $b, Q). $
\ENDFOR
\ENSURE $\widetilde{x}_S$. 
\end{algorithmic}
\end{algorithm}
\begin{algorithm}[t]
\caption{One Stage AccSVRDA $(\widetilde{y}, \widetilde{x}, \eta, m, b, Q)$}
\label{one_Stage_accsvrda}
\begin{algorithmic}
%[1]
\REQUIRE $\widetilde{y}, \widetilde{x}$, $\eta > 0$, $m\in  \mathbb{N}$, $b \in [n]$, $Q$. 
\STATE $x_0 = z_0 = \widetilde{y}$, $\bar{g}_0 = 0$, $\theta_0 = \frac{1}{2}$.
\FOR {$k=1$ to $m$}
\STATE Pick independently $i_k^{1}, \ldots, i_k^{b} \sim Q$, $I_k = \{ i_k^{\ell} \}_{\ell=1}^b$.
\STATE $\theta_k = \frac{k+1}{2}. $
\STATE $y_k = \left(1-\frac{1}{\theta_k}\right)x_{k-1} + \frac{1}{\theta_k}z_{k-1}. $ 
\STATE $g_{k} = \frac{1}{b}\sum_{i \in I_k} \frac{1}{nq_i}\left(\nabla f_{i}(y_k)-\nabla f_{i}(\widetilde{x})\right)+\nabla F(\widetilde{x}). $
\STATE $\bar{g}_k = \left(1-\frac{1}{\theta_k}\right)\bar{g}_{k-1} + \frac{1}{\theta_k}g_k. $ 
\STATE $z_{k}= \underset{z \in \mathbb{R}^d}{\mathrm{argmin}} \left\{ \langle \bar g_k, z \rangle + R(z) + \frac{1}{2\eta \theta_k \theta_{k-1}}\|z-z_0\|^2 \right\}$ 
\STATE \hspace{1.25em}$= \mathrm{prox}_{\eta \theta_k \theta_{k-1}R}\left(z_0-\eta \theta_k \theta_{k-1}\bar{g}_k\right). $ 
\STATE $x_{k}= \left(1-\frac{1}{\theta_k}\right)x_{k-1} + \frac{1}{\theta_k}z_{k}. $
%\label{orda_line}
\ENDFOR
\ENSURE $(x_m, z_m)$.
\end{algorithmic}
\end{algorithm}
\begin{algorithm}[t]
\caption{DASVRDA$^{\mathrm{sc}}(\check x_{0}, \gamma, \{ L_i \}_{i=1}^n , m, b, S, T$)}
\label{dasvrda_sc}
\begin{algorithmic}
%[1]
\REQUIRE $\check x_{0} \in \mathbb{R}^d$, $\gamma \geq 1$, $\{ L_i > 0\}_{i=1}^n$, $m\in  \mathbb{N}$, $b \in [n]$, $S, T\in  \mathbb{N}$. 
\FOR {$t=1$ to $T$}
\STATE $\check x_t =$ DASVRDA$^{\mathrm{ns}}(\check x_{t-1}, \check x_{t-1}, \gamma, \{ L_i \}_{i=1}^n, m, b, S)$. 
\ENDFOR
\ENSURE $\check x_T$.
\end{algorithmic}
\end{algorithm}
\begin{algorithm}[t]
\caption{DASVRDA$^{\mathrm{ns}}$ with warm start $(\widetilde{x}_{0}, \gamma, \{ L_i \}_{i=1}^n, m_0, m, b, U, S$)}
\label{warm_start}
\begin{algorithmic}
\STATE $\widetilde{z}_{0} = \widetilde{x}_{0}$, $\bar L = \frac{1}{n}\sum_{i=1}^n L_i$, $Q = \{ q_i\} = \left\{ \frac{L_i}{n\bar L} \right\}$. 
\FOR{$u=1$ to $U$}
\STATE $m_{u} = \lceil \sqrt{\gamma (m_{u-1}+1)m_{u-1}} \rceil$
\ENDFOR
\STATE $m_U' = \lceil \sqrt{(m_U+1)m_U}/(1-1/\gamma)\rceil$. 
\STATE $\eta = \frac{1}{\left(1 + \frac{\gamma (m_U'+1)}{b}\right)\bar L}$
\FOR{$u=1$ to $U$}
\STATE $(\widetilde{x}_{u}, \widetilde{z}_{u}) = \text{One Stage AccSVRDA}(\widetilde{z}_{u-1}, \widetilde{x}_{u-1}, \eta, m_{u}, $ $b, Q). $
\ENDFOR
\ENSURE DASVRDA$^{\mathrm{ns}}(\widetilde{x}_U, \widetilde{z}_U, \gamma, \{ L_i \}_{i=1}^n, m_U' ,b, S)$. 
\end{algorithmic}
\end{algorithm}
\section{Convergence Analysis of DASVRDA Method}\label{analysis_sec}%%%%%%%%%%%%%%%%%%%%%%%%%%%%%%%%%%%%%%%%%%%%%%%%
In this section, we provide the convergence analysis of our algorithms. 
%Unless otherwise specified, serial computation is assumed. 
First, we consider the DASVRDA$^{\mathrm{ns}}$ algorithm.
\begin{thm}
\label{main_thm_ns}
Suppose that Assumptions \ref{optimal_solution_assump}, \ref{L-smooth_assump} and \ref{proximal_assump} hold. Let $\widetilde{x}_0, \widetilde{z}_0 \in \mathbb{R}^d$, $\gamma \geq 3$, $m \in \mathbb{N}$, $b \in [n]$ and $S \in \mathbb{N}$. Then DASVRDA$^{\mathrm{ns}}(\widetilde{x}_{0}, \widetilde{z}_0, \gamma, \{ L_i \}_{i=1}^n, m, b, S)$ satisfies
\begin{align}
&\mathbb{E} \left[P(\widetilde{x}_{S})-P(x_{*})\right]  \leq \frac{4}{(S+2)^2}\left(P(\widetilde{x}_0)-P(x_*)\right) \notag \\
&+ \frac{8}{\left(1 - \frac{1}{\gamma}\right)^2\eta(S+2)^2(m+1)m}\| \widetilde{z}_0 - x_* \|^2. \notag \notag 
\end{align}
\end{thm}
The proof of Theorem \ref{main_thm_ns} is found in the supplementary material (Section \ref{proof_main_thm_ns_sec}). We can easily see that the optimal choice of $\gamma$ is $(3 + \sqrt{9 + 8b/(m+1)})/2 = O(1+b/m)$ (see Section \ref{proof_main_thm_ns_sec} of supplementary material). We denote this value as $\gamma_*$. Using Theorem \ref{main_thm_ns}, we can establish the convergence rate of DASVRDA$^{\mathrm{ns}}$ with warm start (Algorithm \ref{warm_start}). 
\begin{thm}
\label{main_thm_ns_warm}
Suppose that Assumptions \ref{optimal_solution_assump}, \ref{L-smooth_assump} and \ref{proximal_assump} hold. Let $\widetilde{x}_0 \in \mathbb{R}^d$, $\gamma=\gamma_*$, $m \in \mathbb{N}$, $m_0 = \mathrm{min}\left\{\left\lceil \sqrt{(1+\gamma (m+1)/b)\bar L\frac{\|\widetilde{x}_0 - x_*\|^2}{P(\widetilde{x}_0) - P(x_*)}}\right\rceil, m\right\}$ $\in \mathbb{N}$, $b \in [n]$, $U = \lceil \mathrm{log}_{\sqrt{\gamma}}(m/m_0)\rceil$ and $S \in \mathbb{N}$. Then DASVRDA$^{\mathrm{ns}}$ with warm start$(\widetilde{x}_{0}, \gamma_*, \{ L_i \}_{i=1}^n, $ $m_0, m, b, U, S)$ satisfies
\begin{align}
&\mathbb{E} \left[P(\widetilde{x}_{S})-P(x_{*})\right]  \notag \\
\leq& O\left(\frac{1}{S^2}\left(\frac{1}{m^2} + \frac{1}{mb}\right)\bar L\| \widetilde{x}_0 - x_* \|^2\right). \notag 
\end{align}
\end{thm}
The proof of Theorem \ref{main_thm_ns_warm} is found in the supplementary material (Section \ref{proof_main_thm_ns_warm_sec}). From Theorem \ref{main_thm_ns_warm}, we obtain the following corollary:
%We can choose the optimal value of $\gamma$ based on the following lemma. 
%\newtheorem{lem}{Lemma}
%\begin{lem}\label{optimal_gamma_lemma}
%Define $g(\gamma) \overset{\mathrm{def}}{=} (1 + \gamma(m+1)/b)/(1 - 1/\gamma)^2$ for $\gamma > 1$. Then, 
%$$\gamma_* \overset{\mathrm{def}}{=} \underset{\gamma > 1 }{\mathrm{argmin}} \ g(\gamma) = \frac{1}{2}\left(3 + \sqrt{9+ \frac{8b}{m+1}}\right).$$
%\end{lem}
%Combining Theorem \ref{main_thm_ns} with Lemma \ref{optimal_gamma_lemma} and noting the fact that $\gamma_* \leq O(1+ b/m)$ yield the following corollary:
\begin{cor}\label{main_cor_ns}
Suppose that Assumptions \ref{optimal_solution_assump}, \ref{L-smooth_assump}, and \ref{proximal_assump} hold. Let $\widetilde{x}_0 \in \mathbb{R}^d$,  $\gamma = \gamma_*$ , $m \propto n/b$, $m_0 = \mathrm{min}\left\{\left\lceil \sqrt{(1+\gamma (m+1)/b)\bar L\frac{\|\widetilde{x}_0 - x_*\|^2}{P(\widetilde{x}_0) - P(x_*)}}\right\rceil, m\right\}$ $\in \mathbb{N}$, $b \in [n]$ and $U = \lceil \mathrm{log}_{\sqrt{\gamma}}(m/m_0)\rceil$. If we appropriately choose $S = O(1 +   (1/m + 1/\sqrt{mb})\sqrt{\bar L\| \widetilde{x}_0 - x_* \|^2/\varepsilon})$, then a total computational cost of DASVRDA$^{\mathrm{ns}}$ with warm start$(\widetilde{x}_{0}, \gamma_*, \{ L_i \}_{i=1}^n, m_0, m, b, U, S)$ for $\mathbb{E} \left[P(\widetilde{x}_{S})-P(x_{*})\right] \leq \varepsilon$ is 
{\small $$O\left(d\left( n\mathrm{log}\left(\frac{P(\widetilde{x}_0) - P(x_*)}{\varepsilon}\right) + \left(b + \sqrt{n}\right)\sqrt{\frac{\bar L\| \widetilde{x}_0 - x_* \|^2}{\varepsilon}}\right) \right). $$}
\end{cor}
For the proof of Corollary \ref{main_cor_ns}, see Section \ref{proof_main_cor_ns_sec} of supplementary material.
\begin{rem}
Corollary \ref{main_cor_ns} implies that if the mini-batch size $b$ is $O(\sqrt{n})$,  DASVRDA$^{\mathrm{ns}}$ with warm start$(\widetilde{x}_{0}, \gamma_*, $ $\{ L_i \}_{i=1}^n, m_0, n/b,  b, U, S)$ still achieves the total computational cost of $O( d(n\mathrm{log}(1/\varepsilon) + \sqrt{n\bar L/\varepsilon} ))$, which is better than $O( d(n\mathrm{log}(1/\varepsilon) + \sqrt{nb\bar L/\varepsilon}))$ of Katyusha. 
\end{rem}
\begin{rem}
Corollary \ref{main_cor_ns} also implies that DASVRDA$^{\mathrm{ns}}$ with warm start only needs size $O(\sqrt{n})$ mini-batches for achieving the optimal iteration complexity of  $O(\sqrt{L/\varepsilon})$, when $L/\varepsilon \geq n\mathrm{log}^2(1/\varepsilon)$. In contrast, Katyusha needs size $O(n)$ mini-batches for achieving the optimal iteration complexity. Note that even when $L/\varepsilon \leq n\mathrm{log}^2(1/\varepsilon)$, our method only needs size $\widetilde O(n\sqrt{\varepsilon/L})$ mini-batches\footnotemark\footnotetext{Note that we regard one computation of a full gradient as $n/b$ iterations in size $b$ mini-batch settings.}, that is typically smaller than $O(n)$ of Katyusha.
\end{rem}
Next, we analyze the DASVRDA$^{\mathrm{sc}}$ algorithm for optimally strongly convex objectives. Combining Theorem \ref{main_thm_ns} with the optimal strong convexity of the objective function immediately yields the following theorem, which implies that the DASVRDA$^{\mathrm{sc}}$ algorithm achieves a linear convergence.
\begin{thm}
\label{main_thm_sc}
Suppose that Assumptions \ref{optimal_solution_assump}, \ref{L-smooth_assump}, \ref{proximal_assump} and \ref{osc_assump} hold. Let $\check {x}_0 \in \mathbb{R}^d$, $\gamma = \gamma_*$, $m \in \mathbb{N}$, $b \in [n]$ and $T \in \mathbb{N}$. Define $\rho \overset{\mathrm{def}}{=} $ $4\{ (1-1/\gamma_*)^2 + 4
/(\eta(m+1)m\mu)\}/\{( 1 - 1/\gamma_*)^2(S+2)^2\}$. If $S$ is sufficiently large such that $\rho \in (0, 1)$, then DASVRDA$^{\mathrm{sc}}(\check x_{0}, \gamma_*, \{ L_i \}_{i=1}^n, m, b, S, T)$
satisfies
$$\mathbb{E}[P(\check x_{T}) - P(x_*)] \leq \rho^T [P(\check x_{0}) - P(x_*)].$$
\end{thm}

From Theorem \ref{main_thm_sc}, we have the following corollary.
\begin{cor}\label{main_cor_sc}
Suppose that Assumptions \ref{optimal_solution_assump}, \ref{L-smooth_assump}, \ref{proximal_assump} and \ref{osc_assump} hold. Let $\check {x}_0 \in \mathbb{R}^d$, $\gamma = \gamma_*$, $m \propto n/b$, $b \in [n]$. There exists $S = O(1 + (b/n + 1/\sqrt{n})\sqrt{\bar L/{\mu}})$, such that $1/\mathrm{log}(1/\rho) = O(1)$. Moreover, if we appropriately choose $T = O(\mathrm{log}({P(\check{x}_0) - P(x_*)}/{\varepsilon})$, then a total computational cost of DASVRDA$^{\mathrm{sc}}$$(\check x_{0}, \gamma_*, \{ L_i \}_{i=1}^n, m, b, S, T)$ for $\mathbb{E} \left[P(\check{x}_{T})-P(x_{*})\right] \leq \varepsilon$ is
$$O\left( d\left(n + \left(b + \sqrt{n}\right)\sqrt{\frac{\bar L}{\mu}}\right) \mathrm{log}\left(\frac{P(\check{x}_0) - P(x_*)}{\varepsilon}\right) \right).$$
\end{cor}
\begin{rem}
Corollary \ref{main_cor_sc} implies that if the mini-batch size $b$ is $O(\sqrt{n})$,  DASVRDA$^{\mathrm{sc}}(\check{x}_{0}, \gamma_*, $ $\{ L_i \}_{i=1}^n, n/b, b, S, T)$ still achieves the total computational cost of $O(d (n + \sqrt{n \bar L/\mu}) \mathrm{log}(1/\varepsilon) )$, which is much better than $O( d(n + \sqrt{n b\bar L/\mu}) \mathrm{log}(1/\varepsilon) )$ of APCG, SPDC, and Katyusha. 
\end{rem}
\begin{rem}
Corollary \ref{main_cor_sc} also implies that DASVRDA$^{\mathrm{sc}}$ only needs size $O(\sqrt{n})$ mini-batches for achieving the optimal iteration complexity $O(\sqrt{L/\mu}\mathrm{log}(1/\varepsilon))$, when $L/\mu \geq n$. In contrast, APCG, SPDC and Katyusha need size $O(n)$ mini-batches and AccProxSVRG does $O(\sqrt{L/\mu})$ ones for achieving the optimal iteration complexity. Note that even when $L/\mu \leq n$, our method only needs size $O(n\sqrt{\mu/L})$ mini-batches
\footnote{Note that the required size is $O(n\sqrt{\mu/L}) (\leq O(n))$, which is not
$O(n\sqrt{L/\mu})\geq O(n)$.}. This size is smaller than $O(n)$ of APCG, SPDC, and Katyusha, and the same as that of AccProxSVRG.
\end{rem}

\section{Efficient Implementation for Sparse Data: Lazy Update}\label{lazy_update_sec}
In this section, we briefly comment on the sparse implementations of our algorithms. \par Originally, lazy update was proposed in online settings \cite{duchi2011adaptive}. Generally, it is difficult for accelerated stochastic variance reduction methods to construct lazy update rules because (i) generally, variance reduced gradients are not sparse even if stochastic gradients are sparse; (ii) if we adopt the momentum scheme, the updated solution becomes a convex combination of previous solutions; and (iii) for non-strongly convex objectives, the momentum rate must not be constant. \citet{konevcny2016mini} have tackled the problem of (i) on non-accelerated settings and derived lazy update rules of the ``mini-batch semi-stochastic gradient descent'' (mS2GD) method. \citet{allen2016katyusha} has only mentioned that lazy updates can be applied to Katyusha but did not give explicit lazy update rules of Katyusha. Particularly, for non-strongly convex objectives, it seems to be difficult to derive lazy update rules owing to the difficulty of (iii). The reason we adopt the stochastic dual averaging scheme \cite{xiao2009dual} rather than
 stochastic gradient descent for our method is to be able to overcome the difficulties faced in (i), (ii), and (iii). The lazy update rules of our method support both non-strongly and strongly convex objectives. The formal lazy update algorithms of our method can be found in the supplementary material (Section \ref{supp_lazy_update_sec}). 
\section{Numerical Experiments}\label{experiment_sec}
In this section, we provide numerical experiments to demonstrate the performance of DASVRDA. \par
\begin{figure*}[t]
\begin{subfigmatrix}{3}
\subfigure[a9a, $(\lambda_1, \lambda_2) = (10^{-4}, 0)$]{\includegraphics[width=5.5cm]{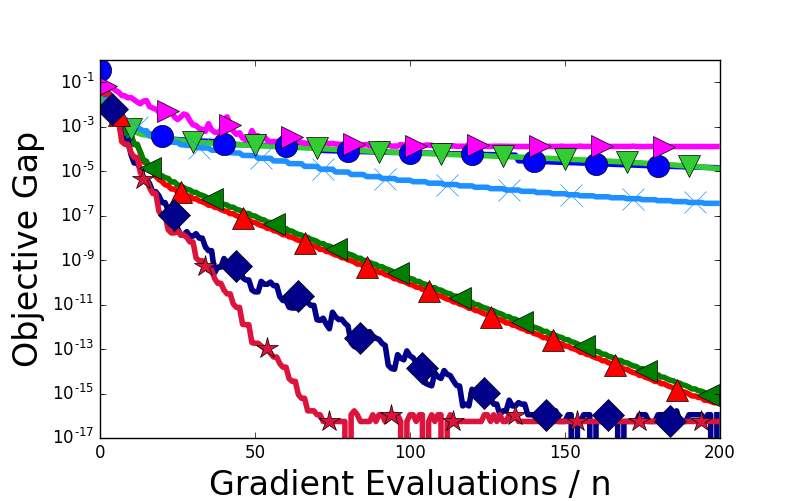}}
\subfigure[a9a, $(\lambda_1, \lambda_2) = (10^{-4}, 10^{-6})$]{\includegraphics[width=5.5cm]{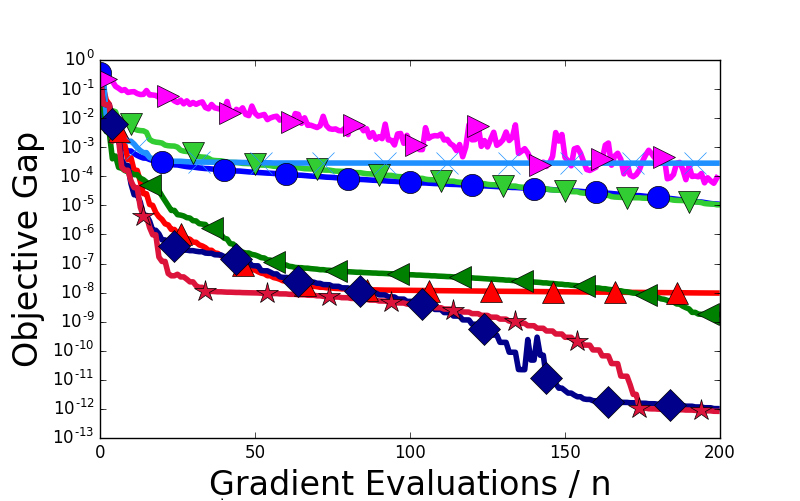}}
\subfigure[a9a, $(\lambda_1, \lambda_2) = (0, 10^{-6})$]{\includegraphics[width=5.5cm]{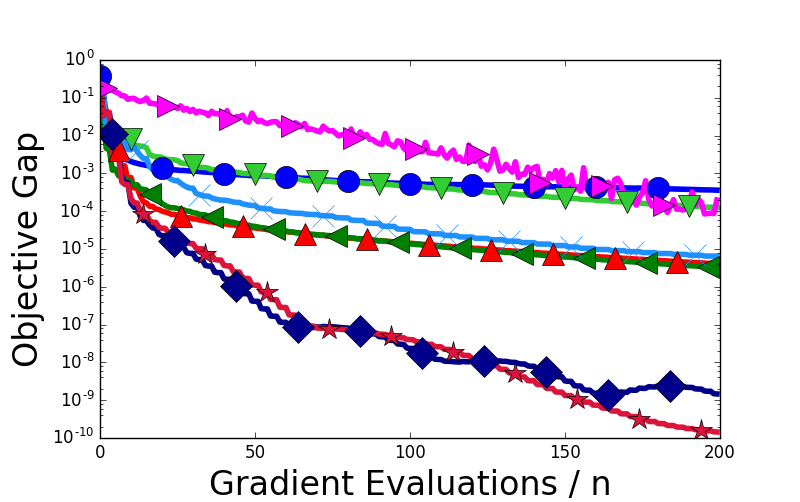}}
\vspace{-0.43cm}
\subfigure[rcv1, $(\lambda_1, \lambda_2) = (10^{-4}, 0)$]{\includegraphics[width=5.5cm]{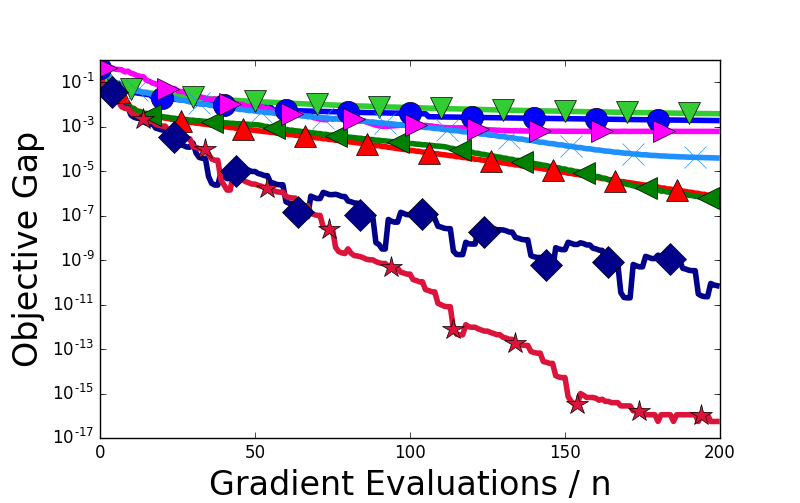}}
\subfigure[rcv1, $(\lambda_1, \lambda_2) = (10^{-4}, 10^{-6})$]{\includegraphics[width=5.5cm]{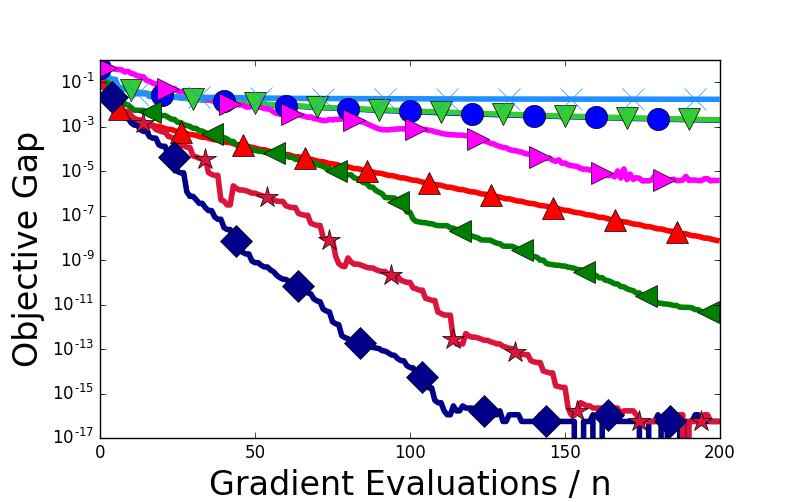}}
\subfigure[rcv1, $(\lambda_1, \lambda_2) = (0, 10^{-6})$]{\includegraphics[width=5.5cm]{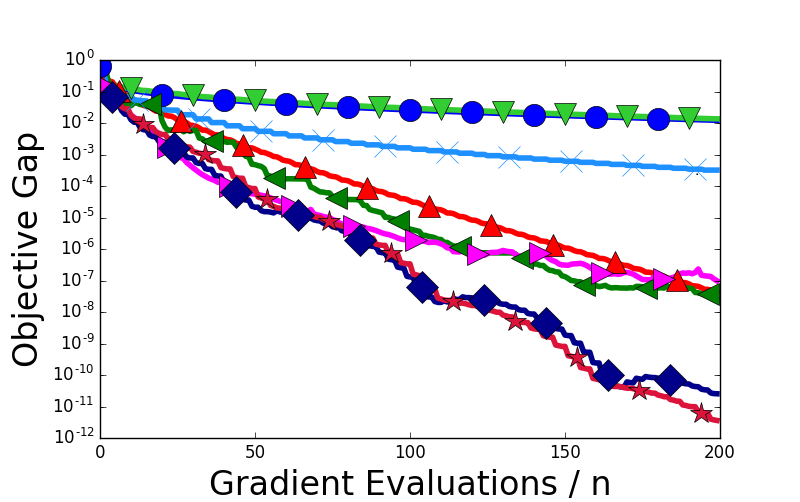}}
\vspace{-0.43cm}
\subfigure[sido0, $(\lambda_1, \lambda_2) = (10^{-4}, 0)$]{\includegraphics[width=5.5cm]{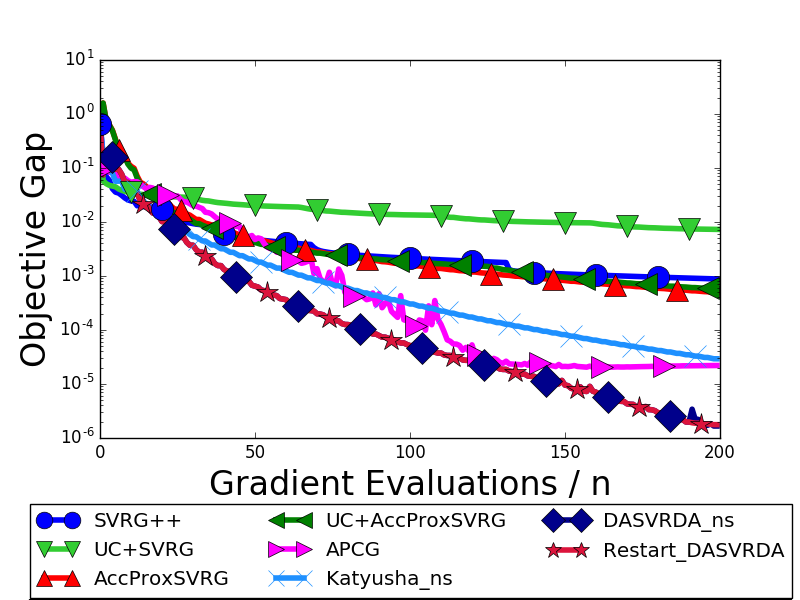}}
\subfigure[sido0, $(\lambda_1, \lambda_2) = (10^{-4}, 10^{-6})$]{\includegraphics[width=5.5cm]{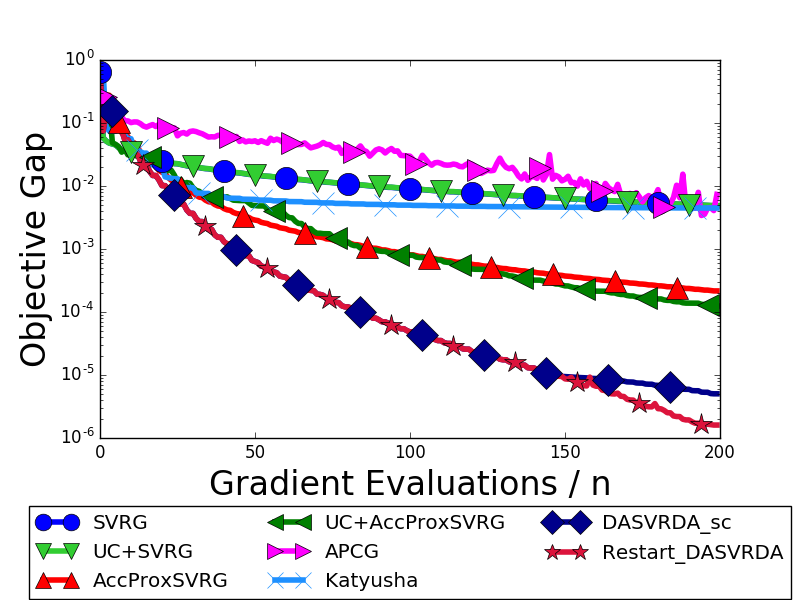}}
\subfigure[sido0, $(\lambda_1, \lambda_2) = (0, 10^{-6})$]{\includegraphics[width=5.5cm]{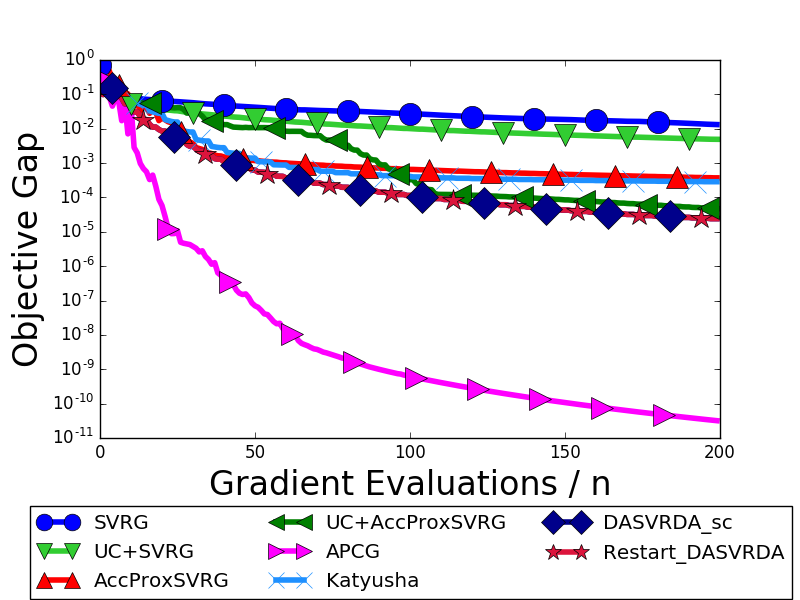}}
\end{subfigmatrix}
\caption{Comparisons on a9a (top), rcv1 (middle) and sido0 (bottom) data sets, for regularization parameters $(\lambda_1, \lambda_2) = (10^{-4}, 0)$ (left), $(\lambda_1, \lambda_2) = (10^{-4}, 10^{-6})$ (middle) and $(\lambda_1, \lambda_2) = (0, 10^{-6})$ (right).}
\label{fig}
\vspace{-0.45cm}
\end{figure*}
We numerically compare our method with several well-known stochastic gradient methods in mini-batch settings: SVRG \cite{xiao2014proximal} (and SVRG$^{++}$ \cite{AllenYang2016}), AccProxSVRG \cite{nitanda2014stochastic}, Universal Catalyst \cite{lin2015universal} , APCG \cite{lin2014accelerated}, and Katyusha \cite{allen2016katyusha}. The details of the implemented algorithms and their parameter tunings are found in the supplementary material. In the experiments, we focus on the regularized logistic regression problem for binary classification, with regularizer $\lambda_1\|\cdot\|_1 + (\lambda_2/2)\|\cdot\|_2^2$. \par
We used three publicly available data sets in the experiments. Their sizes $n$ and dimensions $d$, and common min-batch sizes $b$ for all implemented algorithms are listed in Table \ref{datatable}. \par
\begin{table}[H]
\begin{center}
\begin{tabular}{|c|c|c|c|}\hline
Data sets & $n$ & $d$ & $b$ \\ \hline
a9a  & $32,561$ & $123$ &$180$\\ \hline
rcv1&$20,242$&$47,236$ &$140$\\ \hline
sido0 & $12,678$ & $4,932$ &$100$\\ \hline
\end{tabular}
\end{center}
\caption{Summary of the data sets and mini-batch size used in our numerical experiments}
\label{datatable}
\end{table}
For regularization parameters, we used three settings $(\lambda_1, \lambda_2) = (10^{-4}, 0)$, $(10^{-4}, 10^{-6})$, and $(0, 10^{-6})$. For the former case, the objective is non-strongly convex, and for the latter two cases, the objectives are strongly convex. 
Note that for the latter two cases, the strong convexity of the objectives is $\mu = 10^{-6}$ and is relatively small; thus, it makes acceleration methods beneficial. 
\par
Figure \ref{fig} shows the comparisons of our method with the different methods described above on several settings. ``Objective Gap'' means $P(x)-P(x_*)$ for the output solution $x$. ``Gradient Evaluations $/n$'' is the number of computations of stochastic gradients $\nabla f_i$ divided by $n$. ``Restart$\textunderscore$DASVRDA'' means DASVRDA with heuristic adaptive restarting. We  can observe the following from these results: 
\begin{itemize}
\item Our proposed DASVRDA and Restart$\textunderscore$DASVRDA significantly outperformed all the other methods overall. 
\item DASVRDA with the heuristic adaptive restart scheme efficiently made use of the local strong convexities of non-strongly convex objectives and significantly outperformed vanilla DASVRDA on a9a and rcv1 data sets. For the other settings, the algorithm was still comparable to vanilla DASVRDA. 
\item UC$+$SVRG did not work as well as it did in theory, and its performances were almost the same as that of vanilla SVRG. 
\item UC$+$AccProxSVRG sometimes outperformed vanilla AccProxSVRG but was always outperformed by our methods. 
\item APCG sometimes performed unstably and was outperformed by vanilla SVRG. On sido0 data set, for Elastic Net Setting, APCG significantly outperformed all other methods.  
\item Katyusha outperformed vanilla SVRG overall. However, sometimes Katyusha was slower than vanilla SVRG for Elastic Net Settings. This is probably because SVRG is almost adaptive to local strong convexities of loss functions, whereas Katyusha is not (see the remark in supplementary material). 
%\item DASVRDA had ripples in the trace of the objective gap on a9a and rcv1 data sets. This phenomenon often arises when running determinstic accelerated methods \cite{o2015adaptive}. 
%\item For all data sets, all algorithms for non-strongly convex objectives are much faster than the ones for strongly convex objectives. This is because of the local strongly convexities of the non-strongly convex objectives, particularly the existences of $\ell_1$ regularizers. 
\end{itemize}
\section{Conclusion}
In this paper, we developed a new accelerated stochastic variance reduced gradient method for regularized empirical risk minimization problems in mini-batch settings: DASVRDA. We have shown that DASVRDA achieves the total computational costs of $O(d(n\mathrm{log}(1/\varepsilon) + \sqrt{nL/\varepsilon} + b\sqrt{L/\varepsilon}))$ and $O(d(n+\sqrt{nL/\mu} + b\sqrt{L/\mu})\mathrm{log}(1/\varepsilon))$ in size $b$ mini-batch settings for non-strongly and strongly convex objectives, respectively. In addition, DASVRDA essentially achieves the optimal iteration complexities only with size $O(\sqrt{n})$ mini-batches for both settings.  In the numerical experiments, our method significantly outperformed state-of-the-art methods, including Katyusha and AccProxSVRG. 

\section*{Acknowledgment}
This work was partially supported by MEXT kakenhi (25730013, 25120012,
26280009, 15H01678 and 15H05707), JST-PRESTO and JST-CREST.
\bibliography{bibsvrda}
\bibliographystyle{icml2017}

%\appendix \section{Appendix}

\onecolumn

%\twocolumn[
\icmltitle{Supplementary material: \\ Doubly Accelerated Stochastic Variance Reduced Dual Averaging Method \\for Regularized Empirical Risk Minimization}


\appendix 
In this supplementary material, we give the proofs of Theorem \ref{main_thm_ns}  and the optimality of $\gamma_*$ (Section \ref{proof_main_thm_ns_sec}), Theorem \ref{main_thm_ns_warm} (Section \ref{proof_main_thm_ns_warm_sec}) and Corollary \ref{main_cor_ns} (Section \ref{proof_main_cor_ns_sec}), the lazy update algorithm of our method (Section \ref{supp_lazy_update_sec}) and the experimental details (Section \ref{experimental_sec}). Finally, we briefly discuss DASVRG method, which is a variant of DASVRDA method (Section \ref{dasvrg_sec}).

\section{Proof of Theorem \ref{main_thm_ns}}\label{proof_main_thm_ns_sec}
In this section, we give the comprehensive proof of Theorem \ref{main_thm_ns}. First we analyze One Stage Accelerated SVRDA algorithm. 

\begin{lem}\label{theta_lemma}
The sequence $\{\theta_k \}_{k \geq 1}$ defined in Algorithm \ref{one_Stage_accsvrda} satisfies
$$\theta_k - 1 = \theta_{k-2}$$
for $k \geq 1$, where $\theta_{-1} \overset{\mathrm{def}}{=} 0$.
\end{lem}
\begin{pr}
Since $\theta_k = \frac{k+1}{2}$ for $k \geq 0$, we have that
$$\theta_k -1 = \frac{k+1}{2} - 1 = \frac{k-1}{2} = \theta_{k-2}.$$
\qed
\end{pr}

\begin{lem}\label{theta_product_lemma}
The sequence $\{\theta_k \}_{k \geq 1}$ defined in Algorithm \ref{one_Stage_accsvrda} satisfies
$$\theta_m\theta_{m-1} = \sum_{k=1}^m \theta_{k-1}.$$
\qed
\end{lem}

\begin{pr}
Observe that
$$\theta_m\theta_{m-1} = \frac{m(m+1)}{4} = \sum_{k=1}^m \frac{k}{2} = \sum_{k=1}^m \theta_{k-1}.$$
\qed
\end{pr}

\begin{lem}\label{smoothness_lemma}%%%%%%%%%%%%%%%%%%%%%%%%%%%%%%%%%%%%%%%%%%%%%%
For every $x$, $y \in \mathbb{R}^d$, 
$$F(y) + \langle \nabla F(y), x - y \rangle + R(x)  \leq P(x) - \frac{1}{2\bar L}\frac{1}{n}\sum_{i=1}^{n}\frac{1}{nq_i}\| \nabla f_i(x) - \nabla f_i(y)\|^2.$$
\end{lem}
\begin{pr}
Since $f_i$ is convex and $L_i$-smooth, we have (see \cite{nesterov2013introductory})
$$f_i(y) + \langle \nabla f_i(y), x-y \rangle \leq f_i (x) - \frac{1}{2L_i}\|\nabla f_i(x) - \nabla f_i(y)\|^2. $$
By the definition of $\{ q_i \}$, summing this inequality from $i=1$ to $n$ and dividing it by $n$ results in 
$$F(y) + \langle \nabla F(y), x-y \rangle \leq F(x) - \frac{1}{2\bar L}\frac{1}{n}\sum_{i=1}^{n}\frac{1}{nq_i}\| \nabla f_i(x) - \nabla f_i(y)\|^2.$$
Adding $R(x)$ to the both sides of this inequality gives the desired result. \qed
\end{pr}

\begin{lem}\label{avg_lemma}%%%%%%%%%%%%%%%%%%%%%%%%%%%%%%%%%%%%%%%%%%%%%%%
$$\bar g_k = \frac{1}{\theta_k \theta_{k-1}}\sum_{k' = 1}^{k} \theta_{k'-1} g_{k'} \ (k \geq 1).$$
\end{lem}
\begin{pr}
For $k=1$, $\bar g_1 = g_1 = \frac{1}{1\cdot \frac{1}{2}}\sum_{k' = 1}^1 \frac{1}{2}\cdot g_{k'}$ by the definition of $\theta_0$. \\
Assume that the claim holds for some $k \geq 1$. Then
\begin{align}
\bar g_{k+1} =&\  \left(1-\frac{1}{\theta_{k+1}}\right)\bar g_k + \frac{1}{\theta_{k+1}}g_{k+1} \notag \\
=&\  \left(1-\frac{2}{k+2}\right)\frac{4}{(k+1)k}\sum_{k' = 1}^{k} \theta_{k'-1} g_{k'} + \frac{2}{k+2}g_{k+1} \notag \\
=&\  \frac{4}{(k+2)(k+1)}\sum_{k' = 1}^{k+1} \theta_{k'-1} g_{k'} \notag \\
=&\ \frac{1}{\theta_{k+1}\theta_k}\sum_{k' = 1}^{k+1} \theta_{k'-1} g_{k'}. \notag
\end{align}
The first equality follows from the definition of $\bar g_{k+1}$. Second equality is due to the assumption of the induction.
This finishes the proof for Lemma \ref{avg_lemma}. \qed
\end{pr}

Next we prove the following main lemma for One Stage Accelerated SVRDA. The proof is inspired by the one of AccSDA given in \cite{xiao2009dual}.
\begin{lem}\label{main_lemma}%%%%%%%%%%%%%%%%%%%%%%%%%%%%%%%%%%%%%%%%%%%
Let $\eta < 1/\bar L$. For One Stage Accelerated SVRDA, we have that 
\begin{align}
&\mathbb{E}[P(x_m) - P(x)] \notag \\
\leq&\ \frac{2}{\eta (m+1)m}\| z_{0} - x\|^2 - \frac{2}{\eta(m+1)m}\mathbb{E}\| z_{m} - x\|^2 \notag \\
&+ \frac{2}{(m+1)m}\sum_{k=1}^m \left(\frac{(k+1)k \mathbb{E}\| g_k - \nabla F(y_k)\|^2}{4\left(\frac{1}{\eta}-\bar L\right)} -  \frac{k}{2\bar L}\frac{1}{n}\sum_{i=1}^{n}\frac{1}{nq_i}\mathbb{E}\| \nabla f_i(y_k) - \nabla f_i(x)\|^2 \right), \notag 
\end{align}
for any $x \in \mathbb{R}^d$, where the expectations are taken with respect to $I_k (1 \leq k \leq m)$. 
\end{lem}

\begin{pr}
We define 
\begin{align}
\ell_k (x) &= F(y_{k}) + \langle \nabla F(y_{k}), x-y_{k} \rangle + R(x), \notag \\
\hat \ell_k (x) &= F(y_{k}) + \langle g_k, x-y_{k} \rangle + R(x). \notag
\end{align}

Observe that $\ell_k, \hat \ell_k$ is convex and $\ell_k  \leq P$ by the convexity of $F$ and $R$. Moreover, for $k \geq 1$ we have that
\begin{align}
\sum_{k'=1}^{k} \theta_{k'-1} \hat \ell_{k'} (z) =& \sum_{k' = 1}^{k} \theta_{k'-1} F(y_k) + \sum_{k' = 1}^{k}\langle \theta_{k'-1} g_{k'}, z-y_{k'} \rangle + \sum_{k' = 1}^{k} \theta_{k'-1} R(z) \notag \\
=& \left\langle \theta_k \theta_{k-1} \bar g_k, z \right\rangle + \theta_k \theta_{k-1} R(z) + \sum_{k' = 1}^{k} \theta_{k'-1} F(y_k) - \sum_{k' = 1}^{k} \theta_{k'-1} \langle g_{k'}, y_{k'} \rangle. \notag
\end{align}
The second equality follows from Lemma \ref{avg_lemma} and Lemma \ref{theta_product_lemma}. Thus we see that $z_k = \underset{z \in \mathbb{R}^d} { \mathrm{argmin}\ }\left\{\sum_{k'=1}^{k} \theta_{k'-1} \hat \ell_{k'} (z) + \frac{1}{2\eta}\|z-z_0\|^2 \right\} $. 
Observe that $F$ is convex and $\bar L$-smooth. Thus we have that 
\begin{equation}\label{smoothness_of_F}F(x_k) \leq F(y_{k}) + \langle \nabla F(y_{k}), x_k - y_k \rangle + \frac{\bar L}{2} \|x_k-y_k\|^2. \end{equation}

Hence we see that
\begin{align}
P(x_k) \leq&\  \ell_k(x_k) + \frac{\bar L}{2}\|x_k-y_k\|^2 \notag \\
=&\ \ell_k \left(\left(1-\frac{1}{\theta_k}\right)x_{k-1} + \frac{1}{\theta_k}z_k\right) + \frac{\bar L}{2}\left\| \left(1 - \frac{1}{\theta_k}\right) x_{k-1} + \frac{1}{\theta_k}z_k - y_k\right\|^2 \notag  \\
\leq&\ \left(1 - \frac{1}{\theta_k}\right) \ell_k (x_{k-1}) + \frac{1}{\theta_k} \ell_k(z_k) +\frac{\bar L}{2\theta_k^2}\| z_k - z_{k-1}\|^2 \notag  \\
\leq&\ \left(1 - \frac{1}{\theta_k}\right) P(x_{k-1}) + \frac{1}{\theta_k \theta_{k-1}}\left(\theta_{k-1} \hat \ell_k(z_k) + \frac{\bar L}{2}\| z_k - z_{k-1}\|^2\right)  \notag \\
&- \frac{1}{\theta_k}\langle g_k - \nabla F(y_k), z_k - y_k \rangle \notag \\
=&\ \left(1 - \frac{1}{\theta_k}\right)P(x_{k-1}) + \frac{1}{\theta_k \theta_{k-1}}\left(\theta_{k-1} \hat \ell_k(z_k) + \frac{1}{2\eta}\| z_k - z_{k-1}\|^2\right) \notag \\
&- \frac{1}{2\theta_k \theta_{k-1}} \left(\frac{1}{\eta} - \bar L\right)\| z_k - z_{k-1}\|^2 - \frac{1}{\theta_k}\langle g_k - \nabla F(y_k), z_k - z_{k-1} \rangle \notag \\
&- \frac{1}{\theta_k}\langle g_k - \nabla F(y_k), z_{k-1} - y_k \rangle. \notag 
\end{align}
The first inequality follows from (\ref{smoothness_of_F}). The first equality is due to the definition of $x_k$. The second inequality is due to the convexity of $\ell_k$ and the definition of $y_k$. The third inequality holds because $\ell_k \leq P$ and $\frac{1}{\theta_k^2} \leq \frac{1}{\theta_k \theta_{k-1}}$.

Since $\frac{1}{\eta} > \bar L$, we have that
\begin{align}
 -& \frac{1}{2\theta_k \theta_{k-1}} \left(\frac{1}{\eta} - \bar L\right)\| z_k - z_{k-1}\|^2 - \frac{1}{\theta_k}\langle g_k - \nabla F(y_k), z_k - z_{k-1} \rangle \notag \\ 
\leq& \frac{1}{\theta_k}\frac{\theta_{k-1} \| g_k - \nabla F(y_k)\|^2}{2\left(\frac{1}{\eta}-\bar L\right)} \notag \\
\leq& \frac{\| g_k - \nabla F(y_k)\|^2}{2\left(\frac{1}{\eta}-\bar L\right)}. \notag
\end{align}
The first inequality is due to Young's inequality. The second inequality holds because $\theta_{k-1} \leq \theta_k$.

Using this inequality, we get that
\begin{align}
P(x_k) \leq&\  \left(1 - \frac{1}{\theta_k}\right)P(x_{k-1}) + \frac{1}{\theta_k \theta_{k-1}}\left(\theta_{k-1} \hat \ell_k(z_k) + \frac{1}{2\eta}\| z_k - z_{k-1}\|^2\right) \notag \\
&+ \frac{\| g_k - \nabla F(y_k)\|^2}{2\left(\frac{1}{\eta}-\bar L\right)} - \frac{1}{\theta_k}\langle g_k - \nabla F(y_k), z_{k-1} - y_k \rangle. \notag
\end{align}

Multiplying the both sides of the above inequality by $\theta_k \theta_{k-1}$ yields
\begin{align}
\theta_k \theta_{k-1}P(x_k) \leq&\  \theta_{k-1}(\theta_k -1) P(x_{k-1}) + \theta_{k-1} \hat \ell_k(z_k) + \frac{1}{2\eta}\| z_k - z_{k-1}\|^2 \notag \\
&+ \frac{\theta_k \theta_{k-1}\| g_k - \nabla F(y_k)\|^2}{2\left(\frac{1}{\eta}-\bar L\right)} - \theta_{k-1}\langle g_k - \nabla F(y_k), z_{k-1} - y_k \rangle. \label{ineq_for_dasvrg}
\end{align}

By the fact that $\sum_{k'=1}^{k-1} \theta_{k'-1} \hat \ell_{k'} (z) + \frac{1}{2\eta}\|z-z_0\|^2$ is $\frac{1}{\eta}$-strongly convex and $z_{k-1}$ is the minimizer of $\sum_{k'=1}^{k-1} \theta_{k'-1} \hat \ell_{k'} (z) +\frac{1}{2\eta}\|z-z_0\|^2$ for $k \geq 2$, we have that 
$$\sum_{k'=1}^{k-1} \theta_{k'-1} \hat \ell_{k'} (z_{k-1}) +\frac{1}{2\eta}\|z_{k-1}-z_0\|^2 +  \frac{1}{2\eta}\|z_k -z_{k-1}\|^2 \leq\sum_{k'=1}^{k-1} \theta_{k'-1} \hat \ell_{k'} (z_k) + \frac{1}{2\eta}\|z_k-z_0\|^2$$for $k\geq 1$ (and, for $k=1$, we define $\sum_{k'=1}^0 = 0$). 

Using this inequality, we obtain 
\begin{align}
&\theta_k \theta_{k-1} P(x_k) - \sum_{k'=1}^{k} \theta_{k'-1} \hat \ell_{k'} (z_k) - \frac{1}{2\eta}\|z_k-z_0\|^2 \notag \\
\leq& \theta_{k-1} \theta_{k-2} P(x_{k-1}) -\sum_{k'=1}^{k-1} \theta_{k'-1} \hat \ell_{k'} (z_{k-1}) - \frac{1}{2\eta}\|z_{k-1}-z_0\|^2
+ \frac{\theta_k \theta_{k-1}}{2\left( \frac{1}{\eta} - \bar L\right)}\|g_k- \nabla F(y_k)\|^2 \notag \\ 
&- \theta_{k-1} \langle g_k-\nabla F(y_k), z_{k-1}- y_k \rangle. \notag
\end{align}

Here, the inequality follows from Lemma \ref{theta_lemma} (we defined $\theta_{-1} \overset{\mathrm{def}}{=} 0$).

Summing the above inequality from $k=1$ to $m$ results in
\begin{align}
&\theta_m \theta_{m-1} P(x_{m}) - \sum_{k=1}^{m} \theta_{k-1} \hat \ell_k (z_{m}) - \frac{1}{2\eta}\|z_{m}-z_0\|^2 \notag \\
\leq&\  \sum_{k=1}^{m}\frac{\theta_k \theta_{k-1}\|g_k- \nabla F(y_k)\|^2}{{2\left(\frac{1}{\eta} -  \bar L\right)}} - \sum_{k=1}^{m} \theta_{k-1} \langle g_k-\nabla F(y_k), z_{k-1}- y_k \rangle. \notag
\end{align}

Using $\frac{1}{\eta}$-strongly convexity of  the function $\sum_{k=1}^{m} \theta_{k-1} \hat \ell_k (z) + \frac{1}{2\eta}\|z-z_0\|^2$ and the optimality of $z_{m}$, we have that
$$\sum_{k=1}^{m} \theta_{k-1} \hat \ell_k (z_{m}) + \frac{1}{2\eta}\|z_{m}-z_0\|^2 \leq \sum_{k=1}^{m} \theta_{k-1} \hat \ell_k (x) + \frac{1}{2\eta}\|z_0 - x\|^2 - \frac{1}{2\eta}\|z_{m}-x\|^2.$$

From this inequality, we see that 
\begin{align}
&\theta_m \theta_{m-1} P(x_{m}) \notag \\ \leq&\  \sum_{k=1}^{m} \theta_{k-1} \hat \ell_k (x) + \frac{1}{2\eta}\|z_0-x\|^2 - \frac{1}{2\eta}\|z_{m}-x\|^2 \notag \\
&+\sum_{k=1}^{m}\frac{\theta_k \theta_{k-1}\|g_k- \nabla F(y_k)\|^2}{2\left(\frac{1}{\eta} -  \bar L\right)} - \sum_{k=1}^{m} \theta_{k-1} \langle g_k-\nabla F(y_k), z_{k-1}- y_k \rangle \notag \\
=&\ \sum_{k=1}^{m} \theta_{k-1} \ell_k (x) + \frac{1}{2\eta}\|z_0-x\|^2 - \frac{1}{2\eta}\|z_{m}-x\|^2 \notag \\
&+\sum_{k=1}^{m}\frac{\theta_k \theta_{k-1} \|g_k- \nabla F(y_k)\|^2}{2\left(\frac{1}{\eta} -  \bar L\right)} - \sum_{k=1}^{m} \theta_{k-1} \langle g_k-\nabla F(y_k), z_{k-1}- x \rangle. \notag 
\end{align}

By Lemma \ref{smoothness_lemma} with $x = x$ and $y = y_k$, we have that
$$\ell_k (x) \leq P(x) - \frac{1}{2 \bar L}\frac{1}{n}\sum_{i=1}^{n}\frac{1}{nq_i}\| \nabla f_i(x) - \nabla f_i(y_k)\|^2.$$

Applying this inequality to the above inequality yields
\begin{align}
&\theta_m \theta_{m-1} P(x_{m})  - \sum_{k=1}^m \theta_{k-1} P(x) \notag \\ 
\leq&\ \frac{1}{2\eta}\|z_0-x\|^2 - \frac{1}{2\eta}\|z_{m}-x\|^2 \notag \\
&+ \sum_{k=1}^{m} \left[ \frac{\theta_k \theta_{k-1} \|g_k- \nabla F(y_k)\|^2}{2\left(\frac{1}{\eta} -  \bar L\right)} -\frac{\theta_{k-1}}{2\bar L}\frac{1}{n}\sum_{i=1}^{n}\frac{1}{nq_i}\| \nabla f_i(x) - \nabla f_i(y_k)\|^2\right] \notag \\
&- \sum_{k=1}^{m} \theta_{k-1} \langle g_k-\nabla F(y_k), z_{k-1}- x \rangle. \notag 
\end{align}

Using Lemma \ref{theta_product_lemma} and dividing the both sides of the above inequality by $\theta_m \theta_{m-1}$ result in
\begin{align}
&P(x_{m}) - P(x) \notag \\  
\leq& \frac{1}{2\eta \theta_m \theta_{m-1} }\|z_0-x\|^2 - \frac{1}{2\eta \theta_m \theta_{m-1} }\|z_{m}-x\|^2 \notag \\
&+\frac{1}{\theta_m \theta_{m-1}} \sum_{k=1}^{m} \left[ \frac{\theta_k \theta_{k-1} \|g_k- \nabla F(y_k)\|^2}{2\left(\frac{1}{\eta} -  \bar L\right)} -\frac{\theta_{k-1}}{2\bar L}\frac{1}{n}\sum_{i=1}^{n}\frac{1}{nq_i}\| \nabla f_i(x) - \nabla f_i(y_k)\|^2\right] \notag \\
&- \frac{1}{\theta_m \theta_{m-1}}\sum_{k=1}^{m} \theta_{k-1} \langle g_k-\nabla F(y_k), z_{k-1}- x \rangle. \notag 
\end{align}
Taking the expectations with respect to $I_k (1 \leq k \leq m)$ on the both sides of this inequality yields
\begin{align}
&\mathbb{E}[P(x_{m}) - P(x)] \notag \\  
\leq& \frac{1}{2\eta \theta_m \theta_{m-1} }\|z_0-x\|^2 - \frac{1}{2\eta \theta_m \theta_{m-1} }\mathbb{E}\|z_{m}-x\|^2 \notag \\
&+\frac{1}{\theta_m \theta_{m-1}} \sum_{k=1}^{m} \left[ \frac{\theta_k \theta_{k-1} \mathbb{E}\|g_k- \nabla F(y_k)\|^2}{2\left(\frac{1}{\eta} -  \bar L\right)} -\frac{\theta_{k-1}}{2\bar L}\frac{1}{n}\sum_{i=1}^{n}\frac{1}{nq_i}\mathbb{E} \| \nabla f_i(x) - \nabla f_i(y_k)\|^2\right]. \notag
\end{align}
Here we used the fact that $\mathbb{E}[ g_k-\nabla F(y_k)] = 0$ for $k=1, \ldots , m$. This finishes the proof of Lemma \ref{main_lemma}. 
\qed
\end{pr}

Now we need the following lemma. 
\begin{lem}\label{objective_bound_lemma}%%%%%%%%%%%%%%%%%%%%%%%%%%%%%%%%%%%%%%%
For every $x \in \mathbb{R}^d$, 
$$\frac{1}{n}\sum_{i=1}^{n}\frac{1}{nq_i}\|\nabla f_i(x) - \nabla f_i(x_*)\|^2 \leq 2\bar L (P(x) - P(x_*)).$$
\end{lem}
\begin{pr}
From the argument of the proof of Lemma \ref{smoothness_lemma}, we have
$$\frac{1}{n}\sum_{i=1}^{n}\frac{1}{nq_i}\|\nabla f_i(x) - \nabla f_i(x_*)\|^2 \leq 2\bar L (F(x)-\langle \nabla F(x_*), x-x_* \rangle - F(x_*)).$$
By the optimality  of $x_*$, there exists $\xi_* \in \partial R(x_*)$ such that $\nabla F(x_*) + \xi_* = 0$. Then we have
$$-\langle \nabla F(x_*), x-x_* \rangle = \langle \xi_*, x-x_* \rangle \leq R(x)-R(x_*), $$
and hence 
$$\frac{1}{n}\sum_{i=1}^{n}\frac{1}{nq_i}\| \nabla f_i(x) - \nabla f_i(x_*)\|^2 \leq 2\bar L (P(x) - P(x_*)).$$\qed
\end{pr}

%The following lemma is shown in \cite{nitanda2014stochastic}.
%\begin{lem}\label{variance_lemma}
%Let $\{ \xi_i \}_{i=1}^n$ be a set of vectors in $\mathbb{R}^d$ and $\mu$ denote an average of $\{ \xi_i \}_{i=1}^n$. Let $I$ denote a uniform random variable representing a size $b$ subset of $[n]$. Then, it follows that
%$$\mathbb{E}_{I} \left\| \frac{1}{b}\sum_{i \in I} \xi_i - \mu \right\|^2 = \frac{n-b}{(n-1)b} \mathbb{E}_i \| \xi_i - \mu \|^2.$$
%\end{lem}

\begin{prop}\label{main_prop}
Let $\gamma > 1$ and $\eta \leq 1/((1+\gamma(m+1)/b)\bar L)$. For One Pass Accelerated SVRDA, it follows that 
\begin{align}
\mathbb{E}[P(x_m)-P(\widetilde{x})] \leq \frac{2}{\eta(m+1)m}\mathbb \|\widetilde{y} - \widetilde{x}\|^2 - \frac{2}{\eta(m+1)m}\mathbb{E}\|z_m - \widetilde{x}\|^2,  \notag
\end{align}
and
\begin{align}
&\mathbb{E} [P({x}_{m})-P(x_{*})] \notag \\
\leq&\ \frac{1}{\gamma} (P(\widetilde{x}) - P(x_*)) + \frac{2}{\eta(m+1)m}\| \widetilde{y} - x_*\|^2 - \frac{2}{\eta(m+1)m}\mathbb{E}\|z_m - x_*\|^2, \notag 
\end{align}
where the expectations are taken with respect to $I_k (1 \leq k \leq m)$. 
\end{prop}

\begin{pr}%%%%%%%%%%%%%%%%%%%%%%%%%%%%%%%%%%%%%%%%%%%%%%%%%%%%%%%%%
We bound the variance of the averaged stochastic gradient $\mathbb{E}\|g_k - \nabla F(y_{k})\|^2$: 
\begin{align}
&\mathbb{E}\|g_k - \nabla F(y_{k})\|^2 \notag \\
=&\ \mathbb{E}\left[ \mathbb{E}_{I_k} \|g_k - \nabla F(y_{k})\|^2 \mid [k-1] \right] \notag \\
=&\ \frac{1}{b}\mathbb{E}\left[ \mathbb{E}_{i\sim Q}\| (\nabla f_i(y_k) - \nabla f_i(\widetilde{x}))/nq_i + \nabla F(\widetilde{x}) - \nabla F(y_k) \|^2 \mid [k-1]\right] \notag \\
\leq&\ \frac{1}{b}\mathbb{E}\left[ \mathbb{E}_{i\sim Q}\| (\nabla f_i(y_k) - \nabla f_i(\widetilde{x}))/nq_i\|^2 \mid [k-1]\right] \notag \\
=&\ \frac{1}{b}\mathbb{E}\left[ \frac{1}{n}\sum_{i=1}^n \frac{1}{nq_i} \| \nabla f_i(y_k) - \nabla f_i(\widetilde{x})\|^2\right]  \label{var_bound_by_current_initial} \\
\leq&\  \frac{2}{b}\mathbb{E}\left[ \frac{1}{n}\sum_{i=1}^n \frac{1}{nq_i} \| \nabla f_i(y_k) - \nabla f_i(x_*)\|^2\right] \notag \\
&+   \frac{2}{b}\mathbb{E}\left[ \frac{1}{n}\sum_{i=1}^n \frac{1}{nq_i} \| \nabla f_i(\widetilde{x}) - \nabla f_i(x_*)\|^2\right] \notag \\
\leq&\ \frac{2}{b}\mathbb{E}\left[ \frac{1}{n}\sum_{i=1}^n \frac{1}{nq_i} \| \nabla f_i(y_k) - \nabla f_i(x_*)\|^2\right] +  \frac{4\bar L}{b}(P(\widetilde{x}) - P(x_*)). \label{var_bound_by_optimal_sol} 
\end{align}

The second equality follows from the independency of the random variables $\{ i_1, \ldots, i_b\}$ and the unbiasedness of $(\nabla f_i(y_k) - \nabla f_i(\widetilde{x}))/nq_i + \nabla F(\widetilde{x})$. The first inequality is due to the fact that $\mathbb{E}\|X-\mathbb{E}[X]\|^2 \leq \mathbb{E}\|X\|^2$. The second inequality follows from Young's inequality. The final inequality is due to Lemma \ref{objective_bound_lemma}.

Since $\frac{1}{\eta} \geq \left( 1 + \frac{\gamma(m+1)}{b}\right)\bar L$ and $\gamma > 1$, using (\ref{var_bound_by_current_initial}) yields
\begin{align}
\frac{(k+1)k}{4\left(\frac{1}{\eta} - \bar L\right)}\mathbb{E} \|g_k- \nabla F(y_k)\|^2 -\frac{k}{2\bar L}\mathbb{E} \left[\frac{1}{n}\sum_{i=1}^{n}\frac{1}{nq_i}\| \nabla f_i(y_k) - \nabla f_i(\widetilde{x})\|^2\right]
\leq 0. \notag 
\end{align}

By Lemma \ref{main_lemma} (with $x = \widetilde{x}$) we have 
\begin{align}
\mathbb{E}[P(x_m)-P(\widetilde{x})] \leq \frac{2}{\eta(m+1)m}\mathbb \|\widetilde{y} - \widetilde{x}\|^2 - \frac{2}{\eta(m+1)m}\mathbb{E}\|z_m - \widetilde{x}\|^2. \notag
\end{align}

Similarly, combining Lemma \ref{main_lemma} (with $x=x_*$) with (\ref{var_bound_by_optimal_sol}) results in  
\begin{align}
&\mathbb{E} [P({x}_{m})-P(x_{*})] \notag \\
\leq&\ \frac{1}{\gamma} (P(\widetilde{x}) - P(x_*)) + \frac{2}{\eta(m+1)m}\| \widetilde{y} - x_*\|^2 - \frac{2}{\eta(m+1)m}\mathbb{E}\|z_m - x_*\|^2.  \notag 
\end{align}
These are the desired results.
\qed
\end{pr}

\begin{lem}\label{theta_s_lemma}
The sequence $\{ \widetilde{\theta}_s \}_{s \geq 1}$ defined in Algorithm \ref{dasvrda_ns} satisfies
$$\widetilde{\theta}_{s}\left(\widetilde{\theta}_{s}-1+\frac{1}{\gamma}\right) \leq \widetilde{\theta}_{s-1}^2$$
for any $s \geq 1$.
\end{lem}

\begin{pr}
Since $\widetilde{\theta}_s = \left(1 - \frac{1}{\gamma}\right)\frac{s+2}{2}$ for $s \geq 0$, we have 
\begin{align}
& \widetilde{\theta}_{s}\left(\widetilde{\theta}_{s}-1+\frac{1}{\gamma}\right) \notag \\
=&\ \left( 1 - \frac{1}{\gamma}\right)\frac{s+2}{2}\left( \left(1 - \frac{1}{\gamma}\right)\frac{s+2}{2} - 1 + \frac{1}{\gamma}\right) \notag  \\
=&\ \left( 1 - \frac{1}{\gamma}\right)^2 \frac{s(s+2)}{4} \notag \\
\leq&\ \widetilde{\theta}_{s-1}^2. \notag
\end{align}
This finishes the proof of Lemma \ref{theta_s_lemma}. 
\qed
\end{pr}

Now we are ready to proof Theorem \ref{main_thm_ns}.
\begin{pr_main_thm_ns}
By Proposition \ref{main_prop}, we have
\begin{align}
\mathbb{E}[P(\widetilde{x}_s)-P(\widetilde{x}_{s-1})] \leq \frac{2}{\eta(m+1)m}\mathbb{E} \| \widetilde{y}_{s} - \widetilde{x}_{s-1}\|^2 - \frac{2}{\eta(m+1)m}\mathbb{E}\| \widetilde{z}_{s} - \widetilde{x}_{s-1} \|^2,  \notag
\end{align}
and
\begin{align}
&\mathbb{E} [P(\widetilde{x}_s)-P(x_{*})] \notag \\
\leq&\ \frac{1}{\gamma} \mathbb{E}[P(\widetilde{x}_{s-1}) - P(x_*)] + \frac{2}{\eta(m+1)m}\mathbb{E} \| \widetilde{y}_s - x_*\|^2 - \frac{2}{\eta(m+1)m}\mathbb{E}\| \widetilde{z}_s - x_*\|^2,  \notag
\end{align}
where the expectations are taken with respect to the history of all random variables.

Hence we have
\begin{align}
\mathbb{E}[P(\widetilde{x}_s)-P(\widetilde{x}_{s-1})] \leq \frac{4}{\eta(m+1)m}\mathbb{E} \langle \widetilde{z}_{s} - \widetilde{y}_{s}, \widetilde{x}_{s-1} - \widetilde{y}_{s} \rangle - \frac{2}{\eta(m+1)m}\mathbb{E}\| \tilde{z}_{s} - \tilde{y}_{s} \|^2,  \label{objective_previous_current_diff}
\end{align}
and
\begin{align}
&\mathbb{E} [P(\widetilde{x}_s)-P(x_{*})] \notag \\
\leq&\ \frac{1}{\gamma} \mathbb{E}[P(\widetilde{x}_{s-1}) - P(x_*)] + \frac{4}{\eta(m+1)m}\mathbb{E} \langle \widetilde{z}_s - \widetilde{y}_s, x_* - \widetilde{y}_{s} \rangle - \frac{2}{\eta(m+1)m}\mathbb{E}\| \widetilde{z}_s - \widetilde{y}_s \|^2. \label{objective_error}
\end{align}

Since $\gamma \geq 3$, we have $\widetilde{\theta}_s \geq 1$ for $s \geq 1$. 
Multiplying (\ref{objective_previous_current_diff}) by $\widetilde{\theta}_s(\widetilde{\theta}_s - 1) \geq 0$ and adding  $\widetilde{\theta}_s \times$ (\ref{objective_error}) yield
\begin{align}
&\widetilde{\theta}_s^2 \mathbb{E}[P(\widetilde{x}_s) - P(x_*)] - \widetilde{\theta}_s \left(\widetilde{\theta}_s - 1 + \frac{1}{\gamma}\right)\mathbb{E}[P(\widetilde{x}_{s-1}) - P(x_*)] \notag \\
\leq&\ \frac{4}{\eta(m+1)m}\mathbb{E}\langle \widetilde{\theta}_s(\widetilde{z}_s - \widetilde{y}_s), (\widetilde{\theta}_s -1)\widetilde{x}_{s-1} - \widetilde{\theta}_s \widetilde{y}_{s} + x_* \rangle - \frac{2}{\eta(m+1)m}\mathbb{E}\| \widetilde{\theta}_s(\widetilde{z}_s - \widetilde{y}_s) \|^2. \notag
\end{align}

By Lemma (\ref{theta_s_lemma}), we have 
$$\widetilde{\theta}_s \left(\widetilde{\theta}_s - 1 + \frac{1}{\gamma}\right) \leq \widetilde{\theta}_{s-1}^2$$
for $s \geq 1$.

Thus we get
\begin{align}
&\widetilde{\theta}_s^2 \mathbb{E}[P(\widetilde{x}_s) - P(x_*)] - \widetilde{\theta}_{s-1}^2\mathbb{E}[P(\widetilde{x}_{s-1}) - P(x_*)] \notag \\
\leq&\ \frac{4}{\eta(m+1)m}\mathbb{E}\langle \widetilde{\theta}_s(\widetilde{z}_s - \widetilde{y}_s), (\widetilde{\theta}_s -1)\widetilde{x}_{s-1} - \widetilde{\theta}_s \widetilde{y}_{s} + x_* \rangle - \frac{2}{\eta(m+1)m}\mathbb{E}\| \widetilde{\theta}_s(\widetilde{z}_s - \widetilde{y}_s) \|^2. \notag \\
=&\ \frac{2}{\eta(m+1)m}\left( \mathbb{E}\|(\widetilde{\theta}_s -1)\widetilde{x}_{s-1} - \widetilde{\theta}_s \widetilde{y}_{s} + x_* \|^2 - \mathbb{E}\| (\widetilde{\theta}_s - 1 )\widetilde{x}_{s-1} - \widetilde{\theta}_s\widetilde{z}_s + x_*\|^2\right) \notag 
\end{align}

Since $\widetilde{y}_s = \widetilde{x}_{s-1} + \frac{\widetilde{\theta}_{s-1}-1}{\widetilde{\theta}_s}(\widetilde{x}_{s-1} - \widetilde{x}_{s-2}) + \frac{\widetilde{\theta}_{s-1}}{\widetilde{\theta}_s}(\widetilde{z}_{s-1} - \widetilde{x}_{s-1})$,  
we have
$$(\widetilde{\theta}_s -1)\widetilde{x}_{s-1} - \widetilde{\theta}_s \widetilde{y}_{s} + x_*  = (\widetilde{\theta}_{s-1} - 1)\widetilde{x}_{s-2} - \widetilde{\theta}_{s-1}\widetilde{z}_{s-1} + x_*. $$

Therefore summing the above inequality from $s=1$ to $S$, we obtain 
\begin{align}
&\widetilde{\theta}_s^2\mathbb{E}[P(\widetilde{x}_S) - P(x_*)] \notag \\ 
\leq&\  \widetilde{\theta}_0^2(P(\widetilde{x}_0) - P(x_*)) + \frac{2}{\eta(m+1)m}\| (\widetilde{\theta}_0 - 1)\widetilde{x}_{-1} - \widetilde{\theta}_0 \widetilde{z}_{0} + x_*\|^2 \notag \\
=&\ \left(1-\frac{1}{\gamma}\right)^2(P(\widetilde{x}_0) - P(x_*)) + \frac{2}{\eta(m+1)m}\| \widetilde{z}_{0} - x_*\|^2. \notag
\end{align}

Dividing both sides by $\widetilde{\theta}_s^2$ finishes the proof of Theorem \ref{main_thm_ns}. \qed
\end{pr_main_thm_ns}

\subsection*{Optimal choice of $\gamma$}
We can choose the optimal value of $\gamma$ based on the following lemma.
\begin{lem}\label{optimal_gamma_lemma}
Define $g(\gamma) \overset{\mathrm{def}}{=} \frac{\left(1 + \frac{\gamma(m+1)}{b}\right)}{\left(1 - \frac{1}{\gamma}\right)^2}$ for $\gamma > 1$. Then, 
$$\gamma_* \overset{\mathrm{def}}{=} \underset{\gamma > 1 }{\mathrm{argmin}} \ g(\gamma) = \frac{1}{2}\left(3 + \sqrt{9+ \frac{8b}{m+1}}\right).$$
\end{lem}
\begin{pr}
First observe that
\begin{align}
g'(\gamma) = \frac{ \frac{m+1}{b}\left(1 - \frac{1}{\gamma}\right)^2 - 2\left(1 + \frac{\gamma(m+1)}{b}\right)\left(1 - \frac{1}{\gamma}\right)\frac{1}{\gamma^2} }{\left(1 - \frac{1}{\gamma}\right)^2}. \notag
\end{align}
Hence we have
\begin{align}
&g'(\gamma) = 0 \notag \\
\iff& \frac{m+1}{b}\left(1 - \frac{1}{\gamma}\right)^2 - 2\left(1 + \frac{\gamma(m+1)}{b}\right)\left(1 - \frac{1}{\gamma}\right)\frac{1}{\gamma^2} = 0 \notag \\
\iff& \frac{m+1}{b}(\gamma^2 - \gamma) - 2\left(1 + \frac{\gamma(m+1)}{b}\right) = 0 \notag \\
\iff& \gamma^2 - 3\gamma - \frac{2b}{m+1} = 0 \notag \\
\iff& \gamma = \frac{1}{2}\left(3 + \sqrt{9 + \frac{8b}{m+1}}\right) = \gamma_*. \notag
\end{align}
Here the second and last equivalencies hold from $\gamma > 1$.
Moreover observe that $g'(\gamma) > 0$ for $\gamma > \gamma_*$ and $g'(\gamma) < 0$ for $1 < \gamma < \gamma_*$. This means that $\gamma_*$ is the minimizer of $g$ on the region $\gamma > 1$. 
\qed
\end{pr}

\section{Proof of Theorem \ref{main_thm_ns_warm}}\label{proof_main_thm_ns_warm_sec}
In this section, we give a proof of Theorem \ref{main_thm_ns_warm}.

\begin{pr_main_thm_ns_warm}
Since $\eta = 1/((1+\gamma(m_U'+1)/b)\bar L) \leq 1/((1+\gamma(m_u+1)/b)\bar L)$, from Proposition \ref{main_prop}, we have 
\begin{align}
&\mathbb{E} [P(\widetilde{x}_{u})-P(x_{*})] + \frac{2}{\eta(m_u+1)m_u}\mathbb{E}\|\widetilde{z}_u - x_*\|^2\notag \\
\leq&\ \frac{1}{\gamma} (P(\widetilde{x}_{u-1}) - P(x_*)) + \frac{2}{\eta(m_u+1)m_u}\| \widetilde{z}_{u-1} - x_*\|^2 \notag \\
=&\ \frac{1}{\gamma} \left( P(\widetilde{x}_{u-1}) - P(x_*) + \frac{2\gamma}{\eta(m_u+1)m_u}\| \widetilde{z}_{u-1} - x_*\|^2\right). \notag 
\end{align}
Since $m_u = \lceil \sqrt{\gamma(m_{u-1}+1)m_{u-1}}\rceil$, we have
$$\frac{2\gamma}{\eta(m_u+1)m_u} \leq  \frac{2}{\eta(m_{u-1}+1)m_{u-1}}.$$
Using this inequality, we obtain that
\begin{align}
&\mathbb{E} [P(\widetilde{x}_{U})-P(x_{*})] + \frac{2}{\eta(m_U+1)m_U}\mathbb{E}\|\widetilde{z}_U - x_*\|^2\notag \\
\leq&\ \frac{1}{\gamma} \left(P(\widetilde{x}_{U-1}) - P(x_*) + \frac{2}{\eta(m_{U-1}+1)m_{U-1}}\mathbb{E}\| \widetilde{z}_{U-1} - x_*\|^2\right) \notag \\
\leq&\ \cdots \notag \\
\leq&\ \frac{1}{\gamma^U}\left(P(\widetilde{x}_{0}) - P(x_*) + \frac{2}{\eta(m_{0}+1)m_{0}}\| \widetilde{z}_{0} - x_*\|^2\right) \notag \\
\leq&\ \frac{1}{\gamma^U}\left(P(\widetilde{x}_{0}) - P(x_*) + \frac{2}{\eta(m_{0}+1)m_{0}}\| \widetilde{x}_{0} - x_*\|^2\right) \notag \\
=&\ O\left(\frac{1}{\gamma^U}(P(\widetilde{x}_{0}) - P(x_*)) \right). \notag
\end{align}
The last equality is due to the definitions of $m_0$ and $\eta$, and the fact $m_U' = O(m_U) = O(\sqrt{\gamma}^{U}m_0) = O(m)$ (see the arguments in the proof of Corollary \ref{main_cor_ns}). Since
$$\left(1- \frac{1}{\gamma}\right)^2(m_U'+1)m_U' \geq (m_U+1)m_U,$$
we get
\begin{align}
&\mathbb{E} [P(\widetilde{x}_{U})-P(x_{*})] + \frac{2}{\left(1-\frac{1}{\gamma}\right)^2\eta(m_U'+1)m_U'}\mathbb{E}\|\widetilde{z}_U - x_*\|^2\notag \\
\leq&\ O\left(\frac{1}{\gamma^U}\left(P(\widetilde{x}_{0}) - P(x_*)\right)\right). \notag 
\end{align}
Using the definitions of $U$ and $m_0$ and combining this inequality with Theorem \ref{main_thm_ns}, we obtain that desired result. \qed

\end{pr_main_thm_ns_warm}

\section{Proof of Corollary \ref{main_cor_ns}}\label{proof_main_cor_ns_sec}
In this section, we give a proof of Corollary \ref{main_cor_ns}.

\begin{pr_main_cor_ns}
Observe that the total computational cost at the warm start phase becomes
$$O\left(dnU + db\sum_{u=1}^{U}m_u \right).$$
Since $m_u \leq \sqrt{\gamma}m_{u-1}+\sqrt{\gamma} + 1 \leq \sqrt{\gamma}m_{u-1}+2\sqrt{\gamma} \leq \sqrt{\gamma}^2m_{u-2} + 2\sqrt{\gamma} + 2\sqrt{\gamma}^2\leq \cdots \leq \sqrt{\gamma}^um_0 + 2\sum_{u'=1}^u \sqrt{\gamma}^{u'} = O(\sqrt{\gamma}^um_0)$, we have
$$O\left(dnU + db\sum_{u=1}^{U}m_u \right) = O\left(dnU + db\sqrt{\gamma}^Um_0\right) .$$
Suppose that $m \geq m_0 \sqrt{(P(\widetilde{x}_0) - P(x_*))/\varepsilon}$. Then, this condition implies $U = \lceil \mathrm{log}_{\sqrt{\gamma}}(m/m_0)\rceil \geq \mathrm{log}_{\gamma}((P(\widetilde{x}_0) - P(x_*))/\varepsilon)$. Hence we only need to run $u = O(\mathrm{log}_{\gamma}((P(\widetilde{x}_0) - P(x_*))/\varepsilon)) \leq U$ iterations at the warm start phase and running DASVRDA$^{\mathrm{ns}}$ is not needed. Then the total computational cost becomes 
$$O\left(d\left(n\mathrm{log}\frac{P(\widetilde{x}_0) - P(x_*)}{\varepsilon} + bm_0\sqrt{\frac{P(\widetilde{x_0}) - P(x_*)}{\varepsilon}}\right)\right) \leq O\left(d\left(n\mathrm{log}\frac{P(\widetilde{x}_0) - P(x_*)}{\varepsilon}\right)\right),$$
here we used $mb = O(n)$.
Next, suppose that $m \leq m_0 \sqrt{(P(\widetilde{x}_0) - P(x_*))/\varepsilon}$. In this case, the total computational cost at the warm start phase with full U iterations becomes
$$O\left(d\left(n\mathrm{log}\frac{m}{m_0} + mb\right)\right) \leq O\left(d\left(n\mathrm{log}\frac{P(\widetilde{x}_0) - P(x_*)}{\varepsilon}\right)\right).$$

Finally, using Theorem \ref{main_thm_ns_warm} yields the desired total computational cost. \qed
\end{pr_main_cor_ns}

%\section{Proof of Corollary \ref{main_cor_ns}}
%The first half of the statements is immediate from Theorem 
%First we prove that if $\Omega(\mathrm{min}\{1+L/b\varepsilon, n/b\}) \leq m \leq O(n/b)$, DASVRDA$^{\mathrm{ns}}(\widetilde{x_0}, \gamma_*, \{L_i\}_{i=1}^n, m, b, S)$ achieves a total runtime of 
%$$O\left(d\left(n + (b+\sqrt{n})\sqrt{\frac{\bar L \|\widetilde{x}_0 - x_*\|^2}{\varepsilon}}\right)\right).$$ 

\section{Lazy Update Algorithm of DASVRDA Method}\label{supp_lazy_update_sec}
In this section, we discuss how to efficiently compute the updates of the DASVRDA algorithm for sparse data. Specifically, we derive lazy update rules of One Stage Accelerated SVRDA for the following empirical risk minimization problem:
$$\frac{1}{n}\sum_{i=1}^n \psi_i(a_i^{\top}x) + \lambda_1 \|x\|_1 + \frac{\lambda_2}{2} \|x\|_2^2, \ \ \ \ \lambda_1, \lambda_2 \geq 0$$ 

For the sake of simplicity, we define the one dimensional soft-thresholding operator as follows:
$$\mathrm{soft}(z, \lambda) \overset{\mathrm{def}}{=} \mathrm{sign}\left(z\right) \mathrm{max} \{ \left| z \right| - \lambda, 0 \}, $$
for $z\in \mathbb{R}$. Moreover, in this section, we denote $[z_1, z_2]$ as $\{ z \in \mathbb{Z} \mid z_1 \leq z \leq z_2 \}$ for integers  $z_1, z_2 \in \mathbb{Z}$. The explicit algorithm of the lazy updates for One Stage Accelerated SVRDA  is given by Algorithm \ref{one_Stage_accsvrda_lazy_update}. Let us analyze the iteration cost of the algorithm. Suppose that each feature vector $a_i$ is sparse and the expected number of the nonzero elements is $O(d')$. First note that $|\mathcal{A}_k| = O(bd')$ expectedly if $d' \ll d$. For updating $x_{k-1}$, by Proposition \ref{lazy_update_prop}, we need to compute $\sum_{k' \in K_j^{\pm}} \theta_{k'-2}/(1 + \eta \theta_{k'-1}\theta_{k'-2} \lambda_2)$ and $\sum_{k' \in K_j^{\pm}} \theta_{k'-1}\theta_{k'-2}^2/(1 + \eta \theta_{k'-1}\theta_{k'-2} \lambda_2)$ for each $j \in \mathcal{A}_k$. For this, we first make lists $\{ S_k\}_{k=1}^m = \{\sum_{k'=1}^k \theta_{k' -2}/(1 + \eta \theta_{k'-1}\theta_{k'-2} \lambda_2) \}_{k=1}^m$ and $\{ S_k'\}_{k=1}^m = \{\sum_{k'=1}^k \theta_{k' -1}\theta_{k' -2}^2/(1 + \eta \theta_{k'-1}\theta_{k'-2} \lambda_2) \}_{k=1}^m$ before running the algorithm. This needs only $O(m)$. Note that these lists are not depend on coordinate $j$. Since $K_j^{\pm}$ are sets of continuous integers in $[k_j+2, k]$ or unions of two sets of continuous integers in $[k_j+2, k]$,  we can efficiently compute the above sums. For example, if $K_j^+ = [k_j +2, s_-] \cup [s_+, k]$ for some integers $s_{\pm} \in [k_j+2, k]$, we can compute $\sum_{k' \in K_j^{+}} \theta_{k'-2}/(1 + \eta \theta_{k'-1}\theta_{k'-2} \lambda_2)$ as $S_{s_-} - S_{k_j+1} + S_k - S_{s_+ -1}$ and this costs only $O(1)$. Thus, for computing $x_{k-1}$ and $y_k$, we need only $O(bd')$ computational cost. For computing $g_k$, we need to compute the inner product $a_i^{\top}y_k$ for each $i \in I_k$ and this costs $O(bd')$ expectedly. The expected cost of the rest of the updates is apparently $O(bd')$. Hence, the total expected iteration cost of our algorithm in serial settings becomes $O(bd')$ rather than $O(bd)$. Furthermore, we can extend our algorithm to parallel computing settings. Indeed, if we have $b$ processors, processor $b'$ runs on the set $\mathcal{A}_k^{b'} \overset{\mathrm{def}}{=} \{ j \in [d] \mid a_{i_{b'}, j} \neq 0\}$. Then the total iteration cost per processor becomes ideally $O(d')$. Generally the overlap among the sets $\mathcal{A}_k^{b'}$ may cause latency, however for sufficiently sparse data, this latency is negligible. The following proposition guarantees that Algorithm \ref{one_Stage_accsvrda_lazy_update} is equivalent to Algorithm \ref{one_Stage_accsvrda} when $R(x) = \lambda_1 \|x\|_1 + (\lambda_2/2)\|x\|_2^2$. 
\begin{algorithm}[t]
\caption{Lazy Updates for One Stage AccSVRDA $(\widetilde{y},  \widetilde{x}, \eta, m, b, Q)$}
\label{one_Stage_accsvrda_lazy_update}
\begin{algorithmic}
%[1]
\REQUIRE $\widetilde{y}, \widetilde{x}$, $\eta > 0$, $m\in  \mathbb{N}$, $b \in [n]$, $Q$. 
\STATE $x_0 = z_0 = \widetilde{y}$.
\STATE $g_{0, j}^{\mathrm{sum}} = 0\  (j \in [d])$.
\STATE $\theta_0 = \frac{1}{2}$.
\STATE $k_j = 0\ \ (j \in [d])$. 
\STATE $\widetilde{\nabla} = \nabla F(\widetilde{x})$.
\FOR {$k=1$ to $m$}
\STATE Sample independently $i_1, \ldots, i_b \sim Q$. $I_k = \{ i_1, \ldots, i_b \}$. 
\STATE $\mathcal{A}_k = \{ j \in [d] \mid \exists b' \in [b] :a_{i_{b'}, j} \neq 0 \}. $
\STATE $\theta_k = \frac{k+1}{2}. $
\FOR {$ j \in \mathcal{A}_k$}
\STATE \label{lazy_update_line} 
Update $x_{k-1, j}$, $y_{k, j}$ as in Proposition \ref{lazy_update_prop}.
\ENDFOR
\FOR {$ j \in \mathcal{A}_k$} 
\STATE $g_{k, j}= \frac{1}{b}\sum_{i \in I_t} \frac{1}{nq_i}\left( \psi_{i}'(a_i ^{\top}y_k) a_{i, j}-\psi_{i}' (a_i^{\top}\widetilde{x}) a_{i, j}\right)+\widetilde{\nabla}_j . $
\STATE ${g}_{k, j}^{\mathrm{sum}} = {g}_{k_j, j}^{\mathrm{sum}} + \theta_{k-1} g_{k, j} + \left(\theta_k \theta_{k-1} - \theta_{k_j}\theta_{k_j-1} \right)\widetilde{\nabla}_j. $
\STATE $z_{k, j} = \frac{1}{1 + \eta \theta_k \theta_{k-1}\lambda_2}\mathrm{soft} ( z_{0, j} - \eta g_{k, j}^{\mathrm{sum}},\eta \theta_k \theta_{k-1}\lambda_1). $
\STATE $x_{k, j}= \left(1-\frac{1}{\theta_k}\right)x_{k-1, j} + \frac{1}{\theta_k}z_{k, j}. $
\STATE $k_j =  k. $
\ENDFOR 
\ENDFOR
\ENSURE $(x_m, z_m)$.
\end{algorithmic}
\end{algorithm}

\begin{prop}\label{lazy_update_prop}
Suppose that $R(x) = \lambda_1 \|x\|_1 + \frac{\lambda_2}{2}\|x\|_2^2$ with $\lambda_1, \lambda_2 \geq 0$. Let $j \in [d]$, $k_j \in [m]\cup \{0\}$ and $k \geq k_j+1$. Assume that $\nabla_j f_i(y_{k'}) = \nabla_j f_i(\widetilde{x}) = 0$ for any $i \in [b]$ and $k' \in [k_j+1, k-1]$. In Algorithm \ref{one_Stage_accsvrda}, the following results hold:
\begin{align}
x_{k-1, j}=&\ \begin{cases}
x_{0, j} \hspace{26.8em} (k=1)\\
\frac{\theta_{k_j} \theta_{k_j -1}}{\theta_{k-1} \theta_{k-2}} x_{k_j, j} + \frac{1}{\theta_{k-1}\theta_{k-2}}  \sum_{k' \in K_j^+} \frac{\theta_{k'-2}}{1 + \eta \theta_{k'-1}\theta_{k'-2}\lambda_2}(z_{0, j} - M_{k', j}^+) \hspace{1.6em}(k \geq 2) \\
\hspace{5.9em}+\  \frac{1}{\theta_{k-1}\theta_{k-2}} \sum_{k' \in K_j^-} \frac{\theta_{k'-2}}{1 + \eta \theta_{k'-1}\theta_{k'-2}\lambda_2}(z_{0, j} - M_{k', j}^-)
\end{cases}, \notag \\
y_{k, j} =&\ \begin{cases}
x_{0, j} \hspace{26.8em} (k=1) \\
\left(1 - \frac{1}{\theta_k}\right)x_{k-1, j} +  \frac{1}{\theta_k}\frac{1}{1 + \eta \theta_{k-1} \theta_{k-2}\lambda_2} \times  \ \ \ \ \hspace{11.15em}(k \geq 2)\\ 
\hspace{1em} \mathrm{soft} \left( z_{0, j} - \eta g_{k_j, j}^{\mathrm{sum}} - \eta (\theta_{k-1}\theta_{k-2} - \theta_{k_j}\theta_{k_j -1})\widetilde{\nabla}_j, \eta \theta_{k-1} \theta_{k-2}\lambda_1\right)
\end{cases},  \notag
\end{align}
and 
\begin{align}
z_{k, j} = \frac{1}{1 + \eta \theta_k \theta_{k-1}\lambda_2}\mathrm{soft} ( z_{0, j} - \eta g_{k, j}^{\mathrm{sum}}, \eta \theta_k \theta_{k-1}\lambda_1),  \notag
\end{align}
where 
$$M_{k', j}^{\pm} \overset{\mathrm{def}}{=} \eta\theta_{k' -1}\theta_{k'-2} (\widetilde{\nabla}_j \pm \lambda_1) + \eta g_{k_j, j}^{\mathrm{sum}}- \eta \theta_{k_j}\theta_{k_j -1} \widetilde{\nabla}_j, $$
and $K_j^{\pm} \subset [k_j +2, k]$ are defined as follows: \\
Let $c_1 \overset{\mathrm{def}}{=} \frac{\eta \widetilde{\nabla}_j}{4}$, $c_2 \overset{\mathrm{def}}{=} \frac{\eta \lambda_1}{4}$ and $c_3  \overset{\mathrm{def}}{=} \eta g_{k_j, j}^{\mathrm{sum}}- \eta \theta_{k_j}\theta_{k_j -1} \widetilde{\nabla}_j$ to simplify the notation. Note that $c_2 \geq 0$.
Moreover, we define
\begin{align}
D_{\pm} \overset{\mathrm{def}}{=}&\  (c_1 \pm c_2)^2 + 4(c_1 \pm c_2)(z_{0, j} - c_3),  \notag \\
s_{\pm}^{+} \overset{\mathrm{def}}{=}&\  \frac{c_1 + c_2 \pm \sqrt{D_+}}{c_1 + c_2}, \notag \\
s_{\pm}^{-} \overset{\mathrm{def}}{=}&\  \frac{c_1 - c_2 \pm \sqrt{D_-}}{c_1 - c_2}, \notag
\end{align}
where if $s_{\pm}^{\pm}$ are not well defined, we simply assign $0$ (or any number) to $s_{\pm}^{\pm}$.

$1)$ If $c_1 > c_2$, then 
\begin{align}
K_j^+ \overset{\mathrm{def}}{=}&\ 
\begin{cases}
\emptyset \ \ \ \ \hspace{9.85em} (D_+ \leq 0) \\
[k_j+2, k] \cap [ \lceil s_-^+ \rceil, \lfloor s_+^+ \rfloor ]   \ \ \ \ (D_+ > 0) 
\end{cases}, \notag \\
K_j^- \overset{\mathrm{def}}{=}&\ 
\begin{cases}
[k_j+2, k] \ \ \ \ \hspace{6.15em} (D_- \leq 0) \\
[k_j + 2, \lfloor s_-^- \rfloor ] \cup  [ \lceil s_+^- \rceil, k]   \ \ \ \ (D_- > 0) 
\end{cases}. \notag
\end{align}

$2)$ If $c_1 = c_2$, then　
\begin{align}
K_j^+ \overset{\mathrm{def}}{=}&\ 
\begin{cases}
\emptyset \ \ \ \ \hspace{9.87em}(c_2 = 0 \land z_{0, j} \leq c_3) \\
[k_j+2, k] \ \ \ \  \hspace{6.15em}(c_2 = 0 \land z_{0, j} > c_3) \\
\emptyset \ \ \ \ \hspace{9.85em} (c_2 > 0 \land D_+ \leq 0) \\
[k_j+2, k] \cap [ \lceil s_-^+ \rceil, \lfloor s_+^+ \rfloor ]   \ \ \ \ (c_2 > 0 \land D_+ > 0) 
\end{cases}, \notag \\
K_j^- \overset{\mathrm{def}}{=}&\ 
\begin{cases}
[k_j+2, k] \ \ \ \ \hspace{6.15em} (z_{0, j} < c_3) \\
\emptyset    \ \ \ \ \hspace{9.85em}(z_{0, j} \geq c_3) 
\end{cases}. \notag
\end{align}

$3)$ If $|c_1| < c_2$, then
\begin{align}
K_j^+ \overset{\mathrm{def}}{=}&\ 
\begin{cases}
\emptyset \ \ \ \ \hspace{9.85em} (D_+ \leq 0) \\
[k_j+2, k] \cap [ \lceil s_-^+ \rceil, \lfloor s_+^+ \rfloor ]   \ \ \ \ (D_+ > 0) 
\end{cases}, \notag \\
K_j^- \overset{\mathrm{def}}{=}&\ 
\begin{cases}
\emptyset \ \ \ \ \hspace{9.85em} (D_- \leq 0) \\
[k_j+2, k] \cap [ \lceil s_-^- \rceil, \lfloor s_+^- \rfloor ]    \ \ \ \ (D_- > 0) 
\end{cases}. \notag
\end{align}

$4)$ If $c_1 = -c_2$, then
\begin{align}
K_j^+ \overset{\mathrm{def}}{=}&\ 
\begin{cases}
\emptyset \ \ \ \ \hspace{9.85em}(z_{0, j} \leq c_3) \\
[k_j+2, k] \ \ \ \  \hspace{6.15em}(z_{0, j} > c_3) 
\end{cases} \notag \\
K_j^- \overset{\mathrm{def}}{=}&\ 
\begin{cases}
[k_j+2, k] \ \ \ \ \hspace{6.15em}(c_2 = 0 \land z_{0, j} < c_3) \\
\emptyset  \ \ \ \ \hspace{9.85em}(c_2 = 0 \land z_{0, j} \geq c_3) \\
\emptyset \ \ \ \ \hspace{9.85em}(c_2 >0 \land D_- \leq 0) \\
[k_j+2, k] \cap [ \lceil s_-^- \rceil, \lfloor s_+^- \rfloor ]   \ \ \ \ (c_2 > 0 \land D_- > 0) 
\end{cases}. \notag
\end{align}

$5)$ If $c_1 < -c_2$, then
\begin{align}
K_j^+ \overset{\mathrm{def}}{=}&\ 
\begin{cases}
[k_j +2, k] \ \ \ \ \hspace{6.15em}(D_+ \leq 0) \\
[k_j + 2, \lfloor s_-^+ \rfloor ] \cup [ \lceil s_+^+ \rceil, k] \ \ \ \ (D_+ > 0)
\end{cases},  \notag \\
K_j^- \overset{\mathrm{def}}{=}&\ 
\begin{cases}
\emptyset \ \ \ \ \hspace{9.85em} (D_- \leq 0) \\
[k_j+2, k] \cap [ \lceil s_-^- \rceil, \lfloor s_+^- \rfloor ]   \ \ \ \ (D_- > 0) 
\end{cases}. \notag
\end{align}

\end{prop}

\begin{pr}
First we consider the case $k=1$. Observe that 
$$y_{1, j} = \left(1 - \frac{1}{\theta_1}\right)x_{0, j} + \frac{1}{\theta_1}z_{0, j} = z_{0, j} = x_{0, j},$$
and
\begin{align} z_{1, j} =&\ \frac{1}{1 + \eta \theta_1 \theta_0 \lambda_2} \mathrm{soft}\left(z_{0, j} - \eta \theta_1 \theta_0 \frac{1}{\theta_1}g_{1, j}, \eta \theta_1 \theta_0\lambda_1\right) \notag \\
=&\ \frac{1}{1 + \eta \theta_1 \theta_0 \lambda_2} \mathrm{soft}\left(z_{0, j} - \eta \theta_0 g_{1, j}, \eta \theta_1 \theta_0\lambda_1\right) \notag \\
=&\ \frac{1}{1 + \eta \theta_1 \theta_0 \lambda_2} \mathrm{soft}\left(z_{0, j} - \eta g_{1, j}^{\mathrm{sum}}, \eta \theta_1 \theta_0\lambda_1\right). \notag 
\end{align}
Next we consider the case $k \geq 2$. 
We show that \begin{equation} \label{x_convex_combination_equality}
x_{k-1, j} = \frac{\theta_{k_j} \theta_{k_j -1}}{\theta_{k-1} \theta_{k-2}} x_{k_j, j} + \frac{1}{\theta_{k-1}\theta_{k-2}} \sum_{k' = k_j +2}^k \theta_{k'-2}z_{k'-1}.
\end{equation}
For $k = k_j +1$, (\ref{x_convex_combination_equality}) holds. Assume that (\ref{x_convex_combination_equality}) holds for some $k' \geq k_j + 1$. Then  
\begin{align}
x_{k', j} =&\ \left(1 - \frac{1}{\theta_{k'}}\right)x_{k'-1, j} + \frac{1}{\theta_{k'}}z_{k', j} \notag \\
=&\ \left(1 - \frac{1}{\theta_{k'}}\right) \frac{\theta_{k_j} \theta_{k_j -1}}{\theta_{k'-1} \theta_{k'-2}}x_{k_j, j} + \left(1 - \frac{1}{\theta_{k'}}\right)\frac{1}{\theta_{k'-1}\theta_{k'-2}} \sum_{k'' = k_j +2}^{k'} \theta_{k''-2}z_{k''-1} + \frac{1}{\theta_{k'}}z_{k', j} \notag \\
=&\ \frac{\theta_{k_j} \theta_{k_j -1}}{\theta_{k'} \theta_{k'-1}}x_{k_j, j} + \frac{1}{\theta_{k'}\theta_{k'-1}} \sum_{k'' = k_j +2}^{k'} \theta_{k''-2}z_{k''-1} + \frac{1}{\theta_{k'}}z_{k', j} \notag \\
=&\ \frac{\theta_{k_j} \theta_{k_j -1}}{\theta_{k'} \theta_{k'-1}}x_{k_j, j} + \frac{1}{\theta_{k'}\theta_{k'-1}} \sum_{k'' = k_j +2}^{k'+1} \theta_{k''-2}z_{k''-1}. \notag
\end{align}
The first equality is due to the definition of $x_{k'}$. The second equality follows from the assumption of induction. The third equality holds by Lemma \ref{theta_lemma}.
This shows that (\ref{x_convex_combination_equality}) holds.

Next we show that 
\begin{equation}\label{prox_equality}z_{k'-1, j} = \frac{1}{1 + \eta \theta_{k'-1} \theta_{k'-2}\lambda_2} \mathrm{soft}\left(z_{0, j} - \eta g_{k_j, j}^{\mathrm{sum}} -\eta (\theta_{k' -1}\theta_{k' -2} - \theta_{k_j}\theta_{k_j - 1} )\widetilde{\nabla}_j, \eta \theta_{k'-1} \theta_{k'-2}\lambda_1\right),
\end{equation}
for $k' \in [k_j+2, k]$. 

By the definition of $z_{k'-1}$, we have that 
\begin{align}
z_{k'-1, j} =&\ \mathrm{prox}_{\eta \theta_{k'-1} \theta_{k'-2}R} (z_{0} - \eta \theta_{k'-1} \theta_{k'-2} \bar g_{k'-1})_j \notag \\
=&\ \frac{1}{1 + \eta \theta_{k'-1} \theta_{k'-2}\lambda_2} \mathrm{soft}(z_{0, j} - \eta \theta_{k'-1} \theta_{k'-2} \bar g_{k'-1, j}, \eta \theta_{k'-1} \theta_{k'-2}\lambda_1) \notag
\end{align}

From Lemma \ref{avg_lemma}, we see that
\begin{align}
\theta_{k'-1} \theta_{k'-2} \bar g_{k'-1, j} =&\ \sum_{k'' = 1}^{k'-1} \theta_{k''-1} g_{k'', j} \notag \\
=&\ \sum_{k''=1}^{k_j} \theta_{k''-1} g_{k'', j} + \left(\sum_{k''=k_j+1}^{k'-1} \theta_{k''-1} \right) \widetilde{\nabla}_j\notag \\
=&\ g_{k_j, j}^{\mathrm{sum}} + (\theta_{k'-1}\theta_{k'-2} - \theta_{k_j}\theta_{k_j-1})\widetilde{\nabla}_j. \notag
\end{align}

The first and third equality are due to Lemma \ref{theta_product_lemma}. The second equality holds because $g_{k''-1, j} = \widetilde{\nabla}_j$ for $k'' \in [k_j +1, k-1]$ by the assumption. This shows that (\ref{prox_equality}) holds. Observe that 
\begin{align}
z_{k'-1, j} =&\  \frac{1}{1 + \eta \theta_{k'-1} \theta_{k'-2}\lambda_2} \mathrm{soft}\left(z_{0, j} - \eta g_{k_j, j}^{\mathrm{sum}} -\eta (\theta_{k' -1}\theta_{k' -2} - \theta_{k_j}\theta_{k_j - 1} )\widetilde{\nabla}_j, \eta \theta_{k'-1} \theta_{k'-2}\lambda_1\right) \notag \\
=&\ \frac{1}{1 + \eta \theta_{k'-1} \theta_{k'-2}\lambda_2} \mathrm{sign}\left(z_{0, j} - \eta g_{k_j, j}^{\mathrm{sum}} -\eta (\theta_{k' -1}\theta_{k' -2} - \theta_{k_j}\theta_{k_j - 1} )\widetilde{\nabla}_j\right) \notag \\
&\times \mathrm{max} \left\{ \left| z_{0, j} - \eta g_{k_j, j}^{\mathrm{sum}} -\eta (\theta_{k' -1}\theta_{k' -2} - \theta_{k_j}\theta_{k_j - 1} )\widetilde{\nabla}_j\right| - \eta \theta_{k'-1} \theta_{k'-2}\lambda_1 , 0\right\} \notag \\
=&\ 
\begin{cases}
\frac{1}{1 + \eta \theta_{k'-1} \theta_{k'-2}\lambda_2} (z_{0, j} - M_{k', j}^+) \ \ \ \ (z_{0, j} > M_{k', j}^+) \\
0 \hspace{12.85em}( M_{k', j}^- \leq z_{0, j} \leq M_{k', j}^+ )\\
\frac{1}{1 + \eta \theta_{k'-1} \theta_{k'-2}\lambda_2} (z_{0, j} - M_{k', j}^-) \ \ \ \ (z_{0, j} < M_{k', j}^-)
\end{cases}, \notag
\end{align}
where $M_{k', j}^{\pm} = \eta\theta_{k' -1}\theta_{k'-2} (\widetilde{\nabla}_j \pm \lambda_1) + \eta g_{k_j, j}^{\mathrm{sum}}- \eta \theta_{k_j}\theta_{k_j -1} \widetilde{\nabla}_j$.
We define the real valued functions $M^{\pm}$ as follows:
$$M_j^{\pm}(x) \overset{\mathrm{def}}{=} (c_1\pm c_2)x^2 - (c_1\pm c_2)x + c_3, $$
where $c_1 = \frac{\eta \widetilde{\nabla}_j}{4}$, $c_2 = \frac{\eta \lambda_1}{4}$ and $c_3  = \eta g_{k_j, j}^{\mathrm{sum}}- \eta \theta_{k_j}\theta_{k_j -1} \widetilde{\nabla}_j$
Then we see that $M_j^{\pm}(k') = M_{k', j}^{\pm}$. Let
\begin{align}
D_{\pm} \overset{\mathrm{def}}{=}&\  (c_1 \pm c_2)^2 + 4(c_1 \pm c_2)(z_{0, j} - c_3),  \notag \\
s_{\pm}^{+} \overset{\mathrm{def}}{=}&\  \frac{c_1 + c_2 \pm \sqrt{D_+}}{c_1 + c_2}, \notag \\
s_{\pm}^{-} \overset{\mathrm{def}}{=}&\  \frac{c_1 - c_2 \pm \sqrt{D_-}}{c_1 - c_2}, \notag
\end{align}
where if $s_{\pm}^{\pm}$ are not well defined, we simply assign $0$ (or any number) to $s_{\pm}^{\pm}$.
We can easily show that the following results: \\
$1)$ If $c_1 > c_2$, then 
\begin{align}
z_{0, j} > M_j^+(x)\iff&\ 
\begin{cases}
x \in \emptyset \ \ \ \ \hspace{4.95em} (D_+ \leq 0) \\
s_-^+ < x <  s_+^+    \ \ \ \ \hspace{1.74em}(D_+ > 0) 
\end{cases}, \notag \\
z_{0, j} < M_j^-(x) \iff&\ 
\begin{cases}
x \in \mathbb{R} \ \ \ \ \hspace{4.72em} (D_- \leq 0) \\
x < s_-^- \lor x > s_+^-    \ \ \ \ (D_- > 0) 
\end{cases}. \notag
\end{align}

$2)$ If $c_1 = c_2$, then　
\begin{align}
z_{0, j} > M_j^+(x)\iff&\ 
\begin{cases}
x \in \emptyset \ \ \ \ \hspace{4.95em}(c_2 = 0 \land z_{0, j} \leq c_3) \\
x \in \mathbb{R} \ \ \ \  \hspace{4.72em}(c_2 = 0 \land z_{0, j} > c_3) \\
x \in \emptyset \ \ \ \ \hspace{4.95em} (c_2 > 0 \land D_+ \leq 0) \\
s_-^+ < x < s_+^+    \ \ \ \ \hspace{1.74em} (c_2 > 0 \land D_+ > 0) 
\end{cases}, \notag \\
z_{0, j} < M_j^-(x) \iff&\ 
\begin{cases}
x \in \mathbb{R} \ \ \ \ \hspace{4.72em} (z_{0, j} < c_3) \\
x \in \emptyset    \ \ \ \ \hspace{4.95em}(z_{0, j} \geq c_3) 
\end{cases}. \notag
\end{align}

$3)$ If $|c_1| < c_2$, then
\begin{align}
z_{0, j} > M_j^+(x)\iff&\ 
\begin{cases}
x \in \emptyset \ \ \ \ \hspace{4.95em} (D_+ \leq 0) \\
s_-^+ < x < s_+^+    \ \ \ \ \hspace{1.74em}(D_+ > 0) 
\end{cases}, \notag \\
z_{0, j} < M_j^-(x) \iff&\ 
\begin{cases}
x \in \emptyset \ \ \ \ \hspace{4.95em} (D_- \leq 0) \\
s_-^- < x < s_+^-  \ \ \ \ \hspace{1.74em}(D_- > 0) 
\end{cases}. \notag
\end{align}

$4)$ If $c_1 = -c_2$, then
\begin{align}
z_{0, j} > M_j^+(x)\iff&\ 
\begin{cases}
x \in \emptyset \ \ \ \ \hspace{4.95em}(z_{0, j} \leq c_3) \\
x \in \mathbb{R} \ \ \ \  \hspace{4.72em}(z_{0, j} > c_3) 
\end{cases} \notag \\
z_{0, j} < M_j^-(x) \iff&\ 
\begin{cases}
x \in \mathbb{R} \ \ \ \ \hspace{4.72em}(c_2 = 0 \land z_{0, j} < c_3) \\
x \in \emptyset  \ \ \ \ \hspace{4.95em}(c_2 = 0 \land z_{0, j} \geq c_3) \\
x \in \emptyset \ \ \ \ \hspace{4.95em}(c_2 >0 \land D_- \leq 0) \\
s_-^- < x < s_+^-   \ \ \ \ \hspace{1.74em}(c_2 > 0 \land D_- > 0) 
\end{cases}. \notag
\end{align}

$5)$ If $c_1 < -c_2$, then
\begin{align}
z_{0, j} > M_j^+(x)\iff&\ 
\begin{cases}
x \in \mathbb{R} \ \ \ \ \hspace{4.72em}(D_+ \leq 0) \\
x < s_-^+ \lor x > s_+^+  \ \ \ \ (D_+ > 0)
\end{cases},  \notag \\
z_{0, j} < M_j^-(x) \iff&\ 
\begin{cases}
x \in \emptyset \ \ \ \ \hspace{4.95em} (D_- \leq 0) \\
s_-^- < x < s_+^-   \ \ \ \ \hspace{1.74em}(D_- > 0) 
\end{cases}. \notag
\end{align}

The lazy update rules of $x_{k-1, j}$ is derived by combining (\ref{x_convex_combination_equality}) with these results and noting that $k' \in [k_j + 2, k]$. Finally, combining the definition $y_{k, j} = (1 - 1/\theta_k)x_{k-1, j} + (1/\theta_k)z_{k-1, j}$ with (\ref{prox_equality}) gives the lazy update of $y_{k, j}$. The update rule of $z_{k, j}$ is obvious from the proof of (\ref{prox_equality}). \qed
\end{pr}
\section{Experimental Details}\label{experimental_sec}
In this section, we give the experimental details and also comment on the adaptivity of SVRG to local strong convexity. \par
The details of the implemented algorithms and their parameter tunings were as follows: \par
For non-strongly convex cases (($\lambda_1, \lambda_2) = (10^{-4}, 0)$), 
\begin{itemize}
\item SVRG$^{++}$ \cite{AllenYang2016} with default initial epoch length $m = n/(4b)$ \footnotemark\footnotetext{In \cite{AllenYang2016}, the authors have suggested a default initial epoch length $m = n/4$. Since we used mini-batches with size $b$ in our experiments, it was natural to use $m = n/(4b)$. We made sure that using this epoch length improved the performances in all settings. }. We tuned only the learning rate.  
\item AccProxSVRG \cite{nitanda2014stochastic}. We tuned the epoch length, the constant momentum rate and the learning rate, and additional dummy $\ell_2$ regularizer weight for handling a non-strongly convex objective. 
\item UC \cite{lin2015universal} $+$ SVRG \cite{xiao2014proximal} with default epoch length $m = 2n/b$ \footnotemark\footnotetext{In \cite{xiao2014proximal}, the authors has suggested a default initial epoch length $m = 2n$. Since we used mini-batches with size $b$ in our experiments, it was natural to use $m = 2n/b$. We made sure that using this epoch length improved the performances in all settings. }. We tuned $\kappa$ in \cite{lin2015universal} and the learning rate. We fixed $\eta = 1$ in the algorithm of UC (note that $\eta$ is not learning rate).
\item UC $+$ AccProxSVRG. We tuned $\kappa$ in \cite{lin2015universal}, the epoch length, the constant momentum rate and the learning rate. We fixed $\eta = 1$ in the algorithm of UC (note that $\eta$ is not learning rate).
\item APCG \cite{lin2014accelerated}. We tuned the convexity parameter of the dual objective and the learning rate, and additional dummy $\ell_2$ regularizer weight for handling a non-strongly convex objective.  
\item Katyusha$^{\mathrm{ns}}$ \cite{allen2016katyusha} with default epoch length $m =2n/b$ and Katyusha momentum $\tau_2 = 1/2$ following the suggestion of \cite{allen2016katyusha}. We tuned only the learning rate. We did not adopt AdaptReg scheme because Katyusha with AdaptReg was always a bit slower than vanilla Katyusha in our experiments. 
\item DASVRDA$^{\mathrm{ns}}$ with epoch length $m = n/b$ and $\gamma = \gamma_*$. We tuned only the learning rate.
\item Adaptive Restart DASVRDA with epoch length $m = n/b$ and $\gamma = \gamma_*$. We tuned only the learning rate. We used the gradient scheme for the adaptive restarting, that is we restart DASVRDA$^{\mathrm{ns}}$ if $(\widetilde{y}_s - \widetilde{x}_s)^{\top}(\widetilde{y}_{s+1} - \widetilde{x}_{s}) > 0$.
\end{itemize}

For strongly convex cases ($(\lambda_1, \lambda_2) = (10^{-4}, 10^{-6}), (0, 10^{-6})$), 
\begin{itemize}
\item SVRG \cite{xiao2014proximal} with default epoch length $m = 2n/b$. We tuned only the learning rate.   
\item AccProxSVRG \cite{nitanda2014stochastic}. We tuned the epoch length, the constant momentum rate and the learning rate. 
\item UC \cite{lin2015universal} $+$ SVRG \cite{xiao2014proximal} with default epoch length $m = 2n/b$ \footnotemark[4]. We tuned $\kappa$, $q$ in \cite{lin2015universal}  and the learning rate. 
\item UC $+$ AccProxSVRG. We tuned $\kappa$, $q$ in \cite{lin2015universal}, the epoch length, the constant momentum rate and the learning rate. 
\item APCG \cite{lin2014accelerated}. We tuned the convexity parameter of the dual objective and the learning rate.  
\item Katyusha \cite{allen2016katyusha} with default epoch length $m =2n/b$ and Katyusha momentum $\tau_2 = 1/2$ following the suggestion of \cite{allen2016katyusha}. We tuned $\tau_1$ in \cite{allen2016katyusha} and the learning rate.
\item DASVRDA$^{\mathrm{sc}}$ with epoch length $m = n/b$ and $\gamma = \gamma_*$. We tuned the fixed restart interval $S$ and the learning rate.
\item Adaptive Restart DASVRDA with epoch length $m = n/b$ and $\gamma = \gamma_*$. We tuned only the learning rate. We use the gradient scheme for the adaptive restarting, that is we restart DASVRDA$^{\mathrm{ns}}$ if $(\widetilde{y}_s - \widetilde{x}_s)^{\top}(\widetilde{y}_{s+1} - \widetilde{x}_{s}) > 0$.
\end{itemize}

For tuning the parameters, we chose the values that led to the minimum objective value. 
We selected the learning rates from the set $\{10^{p}, 2 \times 10^{p}, 5 \times 10^{p} \mid p \in \{ 0, \pm 1, \pm 2\} \}$ for each algorithm. We selected the epoch lengths from the set $\{ n \times 10^{-k}, 2n \times 10^{-k}, 5n \times 10^{-k} \mid k \in \{ 0, 1, 2, 3\} \}$ and the momentum rates from the set $\{ 1-10^{-k} \mid k \in \{1, 2, 3, 4\} \}$ for AccProxSVRG. We chose the additional dummy $\ell_2$ regularizer weights from the set $\{ 10^{-k}, 0 \mid k \in \{ 4, 5, 6, 8, 12\} \}$ for AccSVRG and APCG.  We selected $\kappa, q$ from the set $\{ 10^{-k}\mid k \in \{ 1, 2, 3, 4, 5, 6\} \}$ for UC.  We chose the convexity parameter from the set $\{ 10^{-k} \mid k \in \{3, 4, 5, 6, 7\} \}$ for APCG. We selected $\tau_1$ from the set $\{10^{-k}, 2 \times 10^{-k}, 5 \times 10^{-k} \mid k \in \{1, 2, 3\} \}$ for Katyusha. We selected the restart interval from the set  $\{10^{k}, 2 \times 10^{k}, 5 \times 10^{k} \mid k \in \{ 0, 1, 2\} \}$ for DASVRDA$^{\mathrm{sc}}$. \par

We fixed the initial points $0 \in \mathbb{R}^d$ for all algorithms. \par

For a fair comparison, we used uniform sampling for all algorithms, because AccProxSVRG does not support non-uniform sampling. 
\subsection*{Remark on the adaptivity of SVRG to local strong convexity}
For an $L1$-reguralization function, the finite sum is often locally strongly convex (or restricted strongly convex) over the set of sparse solutions\footnotemark (see \cite{agarwal2010fast, agarwal2012stochastic, lin2014adaptive}). For example, the logistic models in our experiments satisfy this property. In our elastic net setting, the regularization parameter of $L1$-regularalization ($10^{-4}$) is much larger than the one of $L2$-reguralization ($10^{-6}$). Thus, exploiting the local strong convexity is important for faster convergence. In \cite{allen2015univr}, the author has proposed a modified SVRG algorithm for strongly convex objectives. The essential difference from vanilla SVRG is only using a weighted average of inner updated solutions for the output solution, rather than using uniform average. Their analysis assumes the strong convexity of the finite sum rather than the regularizer and it is possible that the algorithm exploits the local strong convexities of loss functions. More importantly, their algorithm only uses the strong convexity parameter in the weights for the output solution. Specifically, the output solution has the form $\widetilde{x} = (1/\sum_{k=1}^m (1-\mu\eta)^{-k})\sum_{k=1}^m (1-\mu\eta)^{-k}x_k$, where $\eta = O(1/L)$. Thus if the epoch length $m$ is relatively small, the output solution is almost same as the uniformly averaged one (\citet{allen2015univr} assumes $m = O(L/\mu)$ but we have discovered that using $m=O(n/b)$ still guarantees that the algorithm achieves the convergence rate $O((n+bL/\mu)\mathrm{log}(1/\varepsilon))$, that is same as the one of vanilla SVRG in mini-batch settings). The reason why using $m = O(n/b)$ rather than $m=O(n)$ gave faster convergence in our experiments is probably because using $m=O(n/b)$ ensured that the weighted average was nearer to the uniform average than using $m = O(n)$. Therefore, we can say that vanilla SVRG with epoch length $m = O(n/b)$ is almost adaptive to local strong convexity. In contrast, Katyusha uses the strong convexity parameter in the momentum rate and is quite sensitive. Thus it seems to be difficult for Katyusha to exploit local strong convexity. 
\footnotetext{We say that function $F$ is locally strongly convex with respect to $L1$-regularizer if there exist $\gamma > 0$ and $\tau >0$ such that $F(x) - F(y) - \langle \nabla F(y), x-y\rangle \geq (\gamma/2)\|x - y \|_2^2 - \tau\|x-y\|_1^2$ for any $x, y \in \mathbb{R}^d$. Thus, if $F$ is locally strongly convex, then for any $y \in \mathbb{R}^d$ and $\gamma' < \gamma$, $F$ is $\gamma - \gamma'$-strongly convex over the set $\{x \in \mathbb{R}^d \mid \|x-y\|_1 \leq \sqrt{\gamma'/(2\tau)}\|x -y\|_2\}$. This intuitively means that if $F$ is locally strongly convex, $F$ is strongly convex over the set of sparse solutions whose supports do not fluctuate from each other. }

\section{DASVRG method}\label{dasvrg_sec}
In this section, we briefly discuss a SVRG version of DASVRDA method (we call this algorithm DASVRG) and show that DASVRG has the same rates as DASVRDA. \par In Section \ref{algorithm_sec}, we apply the double acceleration scheme to SVRDA method. We can also apply the one to SVRG. The only difference from DASVRDA is the update of $z_t$ in AccSVRDA (Algorithm \ref{one_Stage_accsvrda}). We take the following update for DASVRG:
\begin{equation}\label{dasvrg_update}
z_{k}= \underset{z \in \mathbb{R}^d}{\mathrm{argmin}} \left\{ \langle g_k, z \rangle + R(z) + \frac{1}{2\eta \theta_{k-1}}\|z-z_{k-1}\|^2 \right\}= \mathrm{prox}_{\eta \theta_{k-1}R}\left(z_{k-1}-\eta \theta_{k-1}{g}_k\right). 
\end{equation}
For the convergence analysis of DASVRG, we only need to show that Lemma \ref{main_lemma} is still valid for this algorithm. 
\begin{pr_main_lem_dasvrg}
From (\ref{ineq_for_dasvrg}) in the proof of Lemma \ref{main_lemma} for DASVRDA, we also have
\begin{align}
\theta_k \theta_{k-1}P(x_k) \leq&\  \theta_{k-1}(\theta_k -1) P(x_{k-1}) + \theta_{k-1} \hat \ell_k(z_k) + \frac{1}{2\eta}\| z_k - z_{k-1}\|^2 \notag \\
&+ \frac{\theta_k \theta_{k-1}\| g_k - \nabla F(y_k)\|^2}{2\left(\frac{1}{\eta}-\bar L\right)} - \theta_{k-1}\langle g_k - \nabla F(y_k), z_{k-1} - y_k \rangle, \notag
\end{align}
because the derivation of this inequality does not depend on the update rule of $z_t$. 
Observe that $z_k = \underset{z \in \mathbb{R}^d}{\mathrm{argmin}} \{ \theta_{k-1} \hat \ell_k(z) + 1/(2\eta)\|z-z_{k-1}\|^2\}$ from (\ref{dasvrg_update}). Since $\theta_{k-1} \hat \ell_k(z) + 1/(2\eta)\|z-z_{k-1}\|^2$ is $\eta$-strongly convex, we have
$$\theta_{k-1} \hat \ell_k(z_k) + \frac{1}{2\eta}\| z_k - z_{k-1}\|^2 + \frac{1}{2\eta}\|z_{k}-x\|^2 \leq \theta_{k-1} \hat \ell_k(x) + \frac{1}{2\eta}\| z_{k-1} - x\|^2.$$
Moreover, using the definitions of $\hat \ell$ and $\ell$, and Lemma \ref{smoothness_lemma}, we have 
\begin{align}
&\hat \ell_k(x) = \ell_k(x) + \langle g_k-\nabla F(y_k), x - y_k \rangle \notag \\
\leq& P(x) - \frac{1}{2 \bar L}\frac{1}{n}\sum_{i=1}^{n}\frac{1}{nq_i}\| \nabla f_i(x) - \nabla f_i(y_k)\|^2 + \langle g_k-\nabla F(y_k), x - y_k \rangle \notag .
\end{align}
Hence, we get
\begin{align}
\theta_k \theta_{k-1}(P(x_k)-P(x)) \leq&\  \theta_{k-1}(\theta_k -1) (P(x_{k-1})-P(x)) + + \frac{1}{2\eta}(\| z_{k-1}-x\|^2 - \|z_{k} - x\|^2) \notag \\
&+ \frac{\theta_k \theta_{k-1}\| g_k - \nabla F(y_k)\|^2}{2\left(\frac{1}{\eta}-\bar L\right)} -\frac{\theta_{k-1}}{2 \bar L}\frac{1}{n}\sum_{i=1}^{n}\frac{1}{nq_i}\| \nabla f_i(x) - \nabla f_i(y_k)\|^2 \notag \\
&- \theta_{k-1}\langle g_k - \nabla F(y_k), z_{k-1} - x \rangle. \notag
\end{align}
Note that $\theta_{k-1}(\theta_{k}-1) \leq \theta_{k-1}\theta_{k-2}$ for $k\geq 2$ and $\theta_1 = 1$.
Finally, summing up the above inequality from $k=1$ to $m$, dividing the both sides by $\theta_m \theta_{m-1}$ and taking expectations with respect to $I_k$ ($1\leq k \leq m$) give the desired result. \qed
\end{pr_main_lem_dasvrg}

%\bibliography{bibsvrda}
%\bibliographystyle{icml2017}

\end{document}